\documentclass[12pt]{article}

\newcommand{\maketitlearxiv}{\maketitle}
\newcommand{\maketitleacm}{}

\usepackage{latexsym}
\usepackage{amsfonts}
\usepackage{amsmath}
\usepackage{amssymb}
\usepackage[colorlinks=true,urlcolor=blue]{hyperref}

\title{Theoretical and Practical Aspects of Space-Time DG-SEM Implementations}

\author
{Lea M. Versbach$^{1\ast}$, Viktor Linders$^{1}$, Robert Klöfkorn$^{1}$, Philipp Birken$^{1}$\\
\\
\normalsize{$^{1}$Centre for Mathematical Sciences, Numerical Analysis, Lund University, Lund, Sweden}\\
\\
\normalsize{$^\ast$lea\_miko.versbach@math.lu.se}
}

\date{}

\usepackage[english]{babel}
\usepackage[utf8]{inputenc}

\usepackage{MnSymbol}
\usepackage{caption}
\usepackage{subcaption}
\usepackage{listings}
\usepackage{tikz}
\usepackage{comment}
\usepackage{enumitem}

\usepackage{graphicx}
\usepackage{subcaption}

\usepackage{tabularx}
\usepackage{caption}
\usepackage{nameref}

\usepackage{algorithm}
\usepackage{algorithmic}
\usepackage{todonotes}
\usepackage{dune}

\usepackage{bbm}

\newcommand{\phispace}{V_h}

\newcommand{\grid}{{\mathcal{T}_h}}
\newcommand{\entity}{E}
\newcommand{\elem}{\entity}
\newcommand{\neig}{K}
\newcommand{\isec}{e}

\newcommand{\RRR}{\mathbb{R}}
\newcommand{\NNN}{{\mathbb{N}}}
\newcommand{\R}{\RRR}

\newcommand{\energy}{\varepsilon}
\newcommand{\pressure}{P}

\newcommand{\Fc}{F_c}
\newcommand{\Fv}{F_v}
\newcommand{\stFc}{\underline{F}_c}
\newcommand{\stFv}{\underline{F}_v}

\newcommand{\fluxF}{{H_c}}
\newcommand{\fluxA}{{H_v}}
\newcommand{\stfluxF}{\underline{H}_c}
\newcommand{\stfluxA}{\underline{H}_v}

\newcommand{\vect}[1]{\boldsymbol{#1}}

\newcommand{\vecv}{\vect{v}}
\newcommand{\vecu}{\vect{u}}

\newcommand{\basefct}{\vect{\psi}}
\newcommand{\tbasefct}{\psi}
\newcommand{\stbasefct}{\underline{\vect{\psi}}}

\newcommand{\sol}{u}
\newcommand{\df}{\vecu_h}
\newcommand{\stdf}{\underline{\vecu}_h}
\newcommand{\ucoeff}{\vecu}
\newcommand{\stucoeff}{\underline{\vecu}}

\newcommand{\Nt}{N_{\tau}}

\newcommand{\oper}[1]{\mathcal{#1}}
\newcommand{\spcoper}{\oper{L}_h}
\newcommand{\stspcoper}{\underline{\oper{L}}_h}

\newcommand{\vjump}[1]{[ \! [ {#1} ] \! ]_{\isec} }

\newcommand{\vaver}[1]{\{ \! \! \{ {#1} \} \! \! \}_{\isec} }

\newcommand{\su}{S(\df)}
\newcommand{\stsu}{S(\stdf)}

\newcommand{\n}{\vect{n}}
\renewcommand{\ne}{\n_{\isec}}

\newcommand{\dual}[1]{\langle \basefct, #1 \rangle}
\newcommand{\stdual}[1]{\langle \stbasefct, #1 \rangle}

\newcommand{\assimulo}{\textsc{Assimulo}\xspace}

\usepackage{amsthm}
\theoremstyle{definition} 
\newtheorem{theorem}{Theorem}[section] 

\newtheorem{definition}[theorem]{Definition}

\newtheorem{remark}[theorem]{Remark}

\newcommand{\dt}{\mathrm{d}\tau}

\renewcommand{\vec}[1]{\mathbf{\boldsymbol{#1}}}
\newcommand{\mat}[1]{\mathbf{\boldsymbol{#1}}}

\usepackage{color}
\definecolor{deepblue}{rgb}{0,0,0.5}
\definecolor{deepred}{rgb}{0.6,0,0}
\definecolor{deepgreen}{rgb}{0,0.5,0}

\definecolor{colKeys}{rgb}{0,0,1} 
\definecolor{colIdentifier}{rgb}{0,0,0} 
\definecolor{colComments}{rgb}{0,1,0.3} 
\definecolor{colString}{rgb}{0,0.5,0} 

\definecolor{dkgreen}{rgb}{0,0.6,0} 
\definecolor{gray}{rgb}{0.5,0.5,0.5} 

\usepackage{listings}
\usepackage{xcolor}

\definecolor{codegreen}{rgb}{0,0.6,0}
\definecolor{codegray}{rgb}{0.5,0.5,0.5}
\definecolor{codepurple}{rgb}{0.58,0,0.82}
\definecolor{backcolour}{rgb}{0.95,0.95,0.92}

\lstdefinestyle{mystyle}{
    backgroundcolor=\color{backcolour},   
    commentstyle=\color{codegreen},
    keywordstyle=\color{magenta},
    numberstyle=\tiny\color{codegray},
    stringstyle=\color{codepurple},
    basicstyle=\ttfamily\footnotesize,
    breakatwhitespace=false,         
    breaklines=true,                 
    captionpos=b,                    
    keepspaces=true,                 
    numbers=left,                    
    numbersep=5pt,                  
    showspaces=false,                
    showstringspaces=false,
    showtabs=false,                  
    tabsize=2
}

\graphicspath{ {./Sections/Figures/} }

\begin{document}

\maketitlearxiv

\begin{abstract}

We discuss two approaches for the formulation and implementation
of space-time discontinuous Galerkin spectral element methods (DG-SEM).
In one, time is treated as an additional coordinate direction and a
Galerkin procedure is applied to the entire problem. In the other,
the method of lines is used with DG-SEM in space and the
fully implicit Runge-Kutta method Lobatto IIIC in time.
The two approaches are mathematically equivalent in the sense
that they lead to the same discrete solution. However, in practice
they differ in several important respects, including the terminology used
to describe them, the structure of the resulting software, and the
interaction with nonlinear solvers. Challenges and merits of the two
approaches are discussed with the goal of providing the practitioner with
sufficient consideration to choose which path to follow. Additionally,
implementations of the two methods are provided as a starting point for
further development. Numerical experiments validate the theoretical accuracy
of these codes and demonstrate their utility, even for 4D problems.

\end{abstract}

\maketitleacm


\section{Introduction}
Typically, partial differential equations are numerically treated with a method of lines ansatz; the spatial directions are discretized first, leaving the time variable continuous. The resulting system of ordinary differential equations is then solved using a numerical method for initial value problems.

An alternative ansatz is to treat the time dimension simply as another coordinate direction, and discretize the whole space-time problem simultaneously, resulting in a fully discrete numerical scheme \cite{Neumuller2013}. This approach has several advantages: Moving boundaries can be treated more easily \cite{VanderVegt2006} and parallelization in time is made possible \cite{Gander2015}. However, it also imposes new challenges since the temporal direction is special and needs to follow a causality principle: The solution at a given time is affected and determined only by the solution at earlier times, never the other way around. An overview of space-time computations in practical engineering applications during the last 25 years can be found in \cite{Tezduyar2019}. 

Here, we consider the discontinuous Galerkin spectral element method (DG-SEM); see e.g. \cite{Birken2021} for an overview and \cite{Kopriva2009} for a detailed exposition. These methods have been very successful for spatial discretizations as they are unstructured, of high order and are very suitable for high performance computing \cite{Krais2021}. Further, DG-SEM fits the so-called Summation-By-Parts Simultaneous-Approximation-Term (SBP-SAT) framework \cite{Carpenter1996,Gassner2013}, implying $L_2$ stability for linear problems. Further, in the last decade, work within this framework has resulted in the development of entropy stable (i.e. nonlinearly stable) discretizations of arbitrarily high order \cite{Fisher2013,Fisher20132}.

Our motivation to consider DG-SEM in a space-time formulation is twofold: Firstly, with a specific choice of numerical fluxes, entropy stability can be extended to incorporate the temporal discretization for hyperbolic conservation laws \cite{Friedrich2019}, thereby resulting in a nonlinearly stable, fully discrete scheme. Secondly, the formulation naturally allows for perfectly scaling parallelization in time, with a speedup equal to the number of discretization points within a time element. There are other approaches for parallelization in time that allow for much larger speedups, but need an initial factor of additional processors before giving any speedup at all \cite{Nievergelt1964}.

There is a strong connection between DG discretizations in time and fully implicit Runge-Kutta (RK) methods: DG-SEM in time using an upwind numerical flux is equivalent to the Lobatto IIIC family of RK methods, in the sense that the two methods give the same numerical solution \cite{Boom2015,Ranocha2019}. This observation lends itself to two very different strategies for implementing DG-SEM in space and time. We can either use the method of lines with DG-SEM in space and Lobatto IIIC in time, or we can use space-time DG-SEM.

While mathematically equivalent, there are important differences between these two approaches: 
\begin{itemize}
\item DG and RK methods have been developed largely independently. Hence, there is a language barrier between these communities, where different terminology is used, e.g. when it comes to order.
\item The two approaches leads to different algebraic systems of linear or nonlinear equation. If solved exactly, these systems have the same solutions. However, in practice these solutions must be approximated, typically using iterative solvers. The interplay between iterative methods and the algebraic systems will in general be different, thus the two methods yield unequal numerical solutions.
\item Implementing the two approaches lead to very different software structure, in particular if we wish to reuse existing software. This implies that various numerical tools and techniques may be more readily accessible in one implementation than the other, depending on whether the DG or the RK approach is chosen.
\end{itemize}

In this paper, we discuss these differences in detail so that practitioners can make an educated choice about which path to follow. Further, we present a code base for the two approaches that may be used as a basis for further development of the methods. In particular, we make use of the open source softwares \dune, the Distributed and Unified Numerics Environment, which is a modular toolbox for solving partial differential equations (PDEs) with grid-based methods \cite{Dedner2020}. We also make use of \assimulo \cite{Andersson2015}, a solver package for initial value problems.

This paper is organized as follows: Following a brief literature review below, our target equation and choice of software is introduced in Section \ref{ch:equations}. In Section~\ref{ch:MOL} we introduce the method of lines approach using DG-SEM with Lobatto IIIC for time stepping. In Section~\ref{ch:STDG} the space-time DG-SEM is described. Throughout, code snippets are included to illustrate the details of the implementations.  Theoretical aspects of the two approaches are discussed in Section~\ref{ch:theory}. In particular, we demonstrate the mathematical equivalence of DG-SEM in time and Lobatto IIIC methods, and relate the terminology employed by the DG and RK communities. Practical aspects of the respective implementations are the subject of Section~\ref{ch:practical}. Here we compare algorithmic and implementation specific requirements and merits of the two approaches. In Section~\ref{ch:experiments} we present numerical experiments that validate some of our discussion points before we finish the article with some concluding remarks in Section~\ref{ch:conclusions}.
The \nameref{Appendix_Installation} contains instructions on how to
install the code discussed in this paper.


\subsection{Further reading}

For the practitioner who wishes to delve deeper into various aspects of the topics discussed in this paper, we here suggest a few places to start for further reading.

The book \cite{Kopriva2009} provides much background material on spatial DG-SEM as well as a guide to its implementation. An overview of entropy stable DG-SEM is given in \cite{Gassner2021} and full mathematical detail is provided in \cite{Carpenter2014}. The theory builds upon the SBP-SAT framework, reviews of which are found in \cite{Fernandez2014,Svard2014}.

For a broad background on implicit Runge-Kutta methods, see \cite{Hairer2010}. An overview of the properties of DG-SEM and other SBP-SAT methods for time integration viewed from the Runge-Kutta perspective is given in \cite{Linders2020}. An evaluation of fully implicit RK methods for use in computational fluid dynamics is given in \cite{Jameson2017}, including discussions of how to solve the nonlinear algebraic systems.

Several nonlinear solvers with application to RK and DG methods have recently been presented in the literature. For Newton-type methods, see e.g. \cite{Pazner2018I,Diosady2017} and the references therein. Solvers utilizing multigrid techniques are discussed and analyzed in \cite{Gander2016,Versbach2021,Franciolini2020}.

The implementation of space-time methods with a focus on the challenge of 4D problems has been studied in \cite{Frontin2021}, and the generation of different 4D space-time meshes have been presented in \cite{Behr2008,Caplan2020}.

\section{Governing Equations and Simulation Software} \label{ch:equations}
\label{sec:equations}
We consider a general class of time dependent nonlinear advection-diffusion-reaction problems
\begin{eqnarray}
\label{eqn:general}
  \partial_t \sol  = \oper{L}(\sol) &:=& - \nabla \cdot\big( \Fc(\sol) -
       \Fv(\sol,\nabla\sol) \big) + S(\sol) \ \ \mbox{ in } \Omega \times (0,T)
\end{eqnarray}
for a vector valued function $\sol\colon\Omega\times(0,T)\to\RRR^r$
with $r\in\mathbb{N}^+$ components. Here, $\Omega \subset \RRR^d$, $d=1,2,3$. Suitable initial and boundary conditions are assumed to be available. $\Fc$ and $\Fv$ describe the convective and viscous fluxes respectively, and $S$ is a source term.
We allow for the possibility that any of the coefficients in the partial differential equation (PDE) \eqref{eqn:general} depend explicitly on the spatial variable $x$ and on time $t$, but to simplify the presentation we suppress this dependency in our notation. 

For the discretization of \eqref{eqn:general} we consider two approaches: The first is a method of lines approach, in which the spatial differential operator is discretized using a DG-SEM approximation, yielding a system of ordinary differential equations (ODEs). This system is then solved using a time stepping scheme. In particular, we consider the Lobatto IIIC family of implicit Runge-Kutta methods.

The second approach is to apply the DG-SEM methodology to the entire equation \eqref{eqn:general}, thereby obtaining a fully implicit DG space-time discretization. 

In the following we will include code snippets to clarify the overall structure of the mathematical formulations at hand and to illustrate how the two approaches can be implemented in an existing code base. We utilize 
\dune \cite{dunereview:20}, which is 
a free and open source software framework for the grid-based numerical solution of
PDEs. 
\dune provides one of the most flexible and comprehensive grid interfaces available, allowing $n$-dimensional grids, which we will use in this paper. 
Additionally, state-of-the-art features such as parallelization, grid
adaptivity and load balancing, and moving grids are supported. 
From the variety of \dune modules available we will make 
use of the Python based front-end for \dunefem \cite{dunefem:10} and
\dunefemdg \cite{dunefemdg:21}, which is able to handle weak forms of PDEs described in
the Unified Form Language (UFL) \cite{ufl:14}. As shown in the next section, the
description of weak forms with UFL is straight forward and easy to use. 
Internally, PDEs described in UFL are translated into C++ code just-in-time, to
ensure that the resulting simulation code is performant. 
For a more detailed description we refer to \cite{dunefem:10, dunefemdg:21} and
the tutorial \cite{femtutorial}. 

The implementation of the Lobatto IIIC method (see \cite{lehsten:21}) 
has been done in \assimulo \cite{Andersson2015}, 
which is also a Python package that can be readily used together
with \dunefem. \assimulo provides a high-level interface  for a wide variety of 
classical and modern solvers of ordinary differential equations. e.g. 
SUNDIALS \cite{Hindmarsh2005} and implicit Runge-Kutta 
solvers \cite{Hairer2009,Hairer2010}. The original codes, which are written in 
FORTRAN, C or Python, are wrapped into \assimulo keeping their original form, 
while the user only needs to interact with the Python interface, where the ODE and the 
initial condition of the problem need to be provided, as well as other additional 
information as e.g. the Jacobian, depending on the solver used.
Existing solver options in \assimulo are for instance explicit and implicit Euler, 
Runge-Kutta34, RADAU5ODE, CVODE, IDA, ODASSL, LSODAR, GLIMDA. 

Comments on how to install \dunefemdg and \assimulo are found in \nameref{Appendix_Installation}.


\section{Method of Lines DG-SEM} \label{ch:MOL}
In this section we describe the method of lines (MOL) approach to discretizing \eqref{eqn:general}. A generic DG method is first presented, followed by the specifications needed to obtain the DG-SEM. Finally, the Lobatto IIIC time stepping method is specified.

\subsection{DG-SEM in Space}

\label{sec:dune}
\label{seq:discretization_spatial} 

Given a tessellation $\grid$ of the computational domain $\Omega$ into elements $E$ with
$\bigcup_{\elem \in \grid} \elem = \Omega$, consider the piecewise
polynomial space
\begin{equation}
\label{eqn:vspace}
   \phispace = \{\vecv\in L^2(\Omega,\RRR^{r}) \; \colon
    \vecv|_{\elem}\in[\mathcal{P}_p(\elem)]^{r}, \ \elem\in\grid\},
      \ \;p \in \NNN,
\end{equation}
where $\mathcal{P}_p(\elem)$ is the space of polynomials whose degree do not exceed $p$.
We let $\Gamma_i$ denote the set of intersections between all pairs of elements in $\grid$ and accordingly $\Gamma$ the set of all
intersections including the boundary of $\Omega$. In \dune, the following commands generate the tessellation $\grid$ and the space $\phispace$: \\

\begin{python}
d = 2 # 1,2,3   
from dune.grid import cartesianDomain, structuredGrid as leafGrid  
# create grid that tessellates $[0,1]^d$ with 10 elements in each coordinate direction
T_h = leafGrid(cartesianDomain([0]*d, [1]*d, [10]*d))  

from dune.fem.space import dglagrangelobatto
p = 3 # polynomial degree
# create DG space with Lagrange basis and Gauss-Lobatto interpolation points
V_h = dglagrangelobatto( T_h, order=p )  
\end{python}

We seek an approximate solution $\df \in \phispace$ by discretizing the spatial operator
$\oper{L}(\sol)$  in \eqref{eqn:general}.
To this end we define for all test functions $\basefct \in \phispace$,
\begin{equation}
\label{convDiscr}
\dual { \spcoper(\df) } := \dual{ K_h(\df) } + \dual{ {I}_h(\df) }.
\end{equation}
Here, the element integrals are given by
\begin{eqnarray}
\label{eqn:elementint}
   \dual{ {K}_h(\df) } &:=&
      \sum_{\elem \in \grid} \int_{\elem}
      \big( ( \Fc(\df) - \Fv(\df, \nabla \df ) ) : \nabla\basefct + \su
      \cdot \basefct \big) \, dx,
\end{eqnarray}
where $ : $ denotes the inner product of two second order tensors.
In the code this looks as follows: \\

\begin{python}
from ufl import TrialFunction, TestFunction, inner, grad, dx
# trial and test function
u   = TrialFunction(V_h)
psi = TestFunction(V_h)
# element integral from equation \eqref{eqn:elementint}
K_h = inner(F_c(u) - F_v(u)*grad(u)), grad(psi)) * dx \  # fluxes
    + inner(S(u), psi) * dx                              # source term
\end{python}

The surface
integrals are given by
\begin{eqnarray}
\label{eqn:surfaceint}
   \dual{ {I}_h(\df) } &:=&
      \sum_{e \in \Gamma_i} \int_e \big(
      \vaver{\Fv(\df, \vjump{\df} )^T : \nabla\basefct} +
      \vaver{\Fv(\df, \nabla\df)} : \vjump{\basefct} \big) \, dS \nonumber \\
    &&- \sum_{e \in \Gamma} \int_e \big( \fluxF(\df) - \fluxA(\df,\nabla\df)\big) :
      \vjump{\basefct} \, dS.
\end{eqnarray}

This formulation arises from considering the weak form of the problem: Replace $\sol$ by $\df$ in \eqref{eqn:general}, 
multiply by the test function $\basefct$ and integrate the spatial terms by parts. 
Here, $\fluxF$ and $\fluxA$ are suitable numerical fluxes, imposed at the element interface $e$. 
Further, $\vaver{ \sol }$ and $\vjump{ \sol }$ denote the average and jump of $\sol$ over $e$,  
\begin{equation}
  \vaver{ \sol } := \frac{1}{2} ( \sol_{\elem} + \sol_{\neig} ) \quad \mbox{ and } \quad 
  \vjump{ \sol } := \ne \cdot ( \sol_{\elem} - \sol_{\neig} ) 
\end{equation}
where $\elem$ and $\neig$ are neighboring elements over intersection $\isec$ and
$\ne$ is outward pointing from element $\elem$. 

The corresponding code reads: \\

\begin{python}
from ufl import FacetNormal, FacetArea, CellVolume, avg, jump, dS, ds
# normal and mesh width
n = FacetNormal(V_h)
h_e = avg( CellVolume(V_h) ) / FacetArea(V_h)
# surface integral from equation \eqref{eqn:surfaceint}
I_h = inner(jump(H_c(u), jump(psi)) * dS \  # interior skeleton for convective part
    + H_cb(u)*psi*ds \                        # domain boundary for convective part
    - inner(jump(F_v(u),n),avg(grad(psi)))   * dS \   # symmetry term
    - inner(avg(F_v(u)*grad(u)),jump(psi,n)) * dS \   # consistency term
    + eta/h_e*inner(jump(u, avg(F_v(u))*n),jump(psi,n)) * dS  # penalty term
\end{python}


To obtain the DG-SEM we follow \cite{kopriva:02,kopriva:10}. First, we restrict our focus to cuboid meshes and map each $\elem\in\grid$ to a reference element using an affine mapping. 
In the \dune implementation, the reference element is $[0,1]^d$. This is due to
a generic construction of reference elements of different shapes in arbitrary
dimensions in \dune; see \cite{dedner:12} for details.

In each spatial dimension, a set of $p+1$ Legendre-Gauss-Lobatto (LGL) nodes are introduced and a corresponding set of Lagrange basis polynomials are defined. The discrete solution $\df(t) \in \phispace$ takes the form 
\[\df(t,x) = \sum_i \sol_i(t)\basefct_i(x),\] 
where the sum is taken over all tensor product LGL nodes in $d$ dimensions and $\basefct_i(x)$ is constructed as the product of Lagrange basis polynomials along each dimension. 
This is achieved through the command \\

\begin{python}
# create discrete function given a discrete space  
u_h = V_h.function(name="u_h") 
\end{python}

The convective and viscous fluxes are approximated using the interpolation 
\[\vec{F}_h(t,x) \approx \sum_{i=1} F(u_i(t)) \basefct_i(x),\]
where $F$ is either $F_c$ or $F_v$. A variety of implementations for $F_c$ and $F_v$ is provided by the
\code{dolfin_dg} package (see \cite{Houston:18}), which we use for the Euler equations. 

Finally, the element and surface integrals in \eqref{eqn:elementint} and \eqref{eqn:surfaceint} are approximated using Gauss-Lobatto quadrature rules. The collocation of the quadrature with the LGL nodes results in a diagonal positive definite local mass matrix. The choice of a cuboid mesh and a tensor product formulation of the basis functions ensures that the global mass matrix remains diagonal positive definite and is consequently trivially invertible.

The convective numerical flux $\fluxF$ can be any appropriate numerical flux known for standard finite volume methods. We use the local Lax-Friedrichs (Rusanov) flux function
\begin{equation}
  \label{flux:llf}
  \fluxF^{\! \! LLF}(\df)|_{\isec} :=
      \vaver{ \Fc(\df) } + \frac{\lambda_{\isec}}{2}
         \vjump{ \df }
\end{equation}
where $\lambda_{\isec}$ is an estimate of the maximum wave speed on the interface $\isec$.
Other options are implemented in \dunefemdg (cf.
\cite{dunefemdg:17,dunefemdg:21}) as well.

A wide range of diffusion fluxes $\fluxA$ can be found in the
literature (cf. \cite{cdg2:12} and references therein), 
however, only fluxes from the Interior Penalty family can currently 
be described in UFL due to the missing description and 
implementation in UFL of lifting terms needed for the other fluxes.
For the Interior Penalty method the flux is chosen to be 
\begin{equation}
  \fluxA^{\!\! IP}(\sol, \nabla \sol) = \vaver{\nabla \sol } - \frac{\eta}{h_{\isec}}\vaver{\Fv(\sol, \nabla \sol)} \vjump{\sol}
\end{equation} 
with $\eta$ being the penalty parameter.

\subsection{Temporal Discretization}
\label{TimeDisc}

After spatial discretization, we obtain a system of ODEs for the coefficient functions $\ucoeff(t) = (\sol_1(t), \sol_2(t), \dots)^\top$, which reads
\begin{eqnarray}
  \label{eqn:ode}
  \ucoeff'(t) = \vec{F}(t,\ucoeff(t)), \quad t \in (0,T], \quad \ucoeff(0) = \vec{u}_0.
\end{eqnarray}
Here, $\vec{F}(t,\ucoeff(t)) = \mat{M}^{-1}\spcoper(\df(t))$, where $\spcoper$ is defined in \eqref{convDiscr} and $\mat{M}$ is the (diagonal) global mass matrix of the DG-SEM discretization. The initial data $\vec{u}_0$ for \eqref{eqn:ode} is given by the projection of $\sol_0$ onto $\phispace$.

Any Runge-Kutta method 
can in principle be used to solve \eqref{eqn:ode}. Explicit methods are easy to implement
but suffer from severe time step restrictions for stiff systems. 

Consider instead an implicit RK method with Butcher tableau
$$
\renewcommand\arraystretch{1.2}
\begin{array}
{c|c}
\vec{c} & \mat{A} \\
\hline
& \vec{b}^\top 
\end{array}
$$
The \emph{stage equations} of the RK method take the form
\begin{equation} \label{eq:RK_stage_equations}
\underline{\ucoeff} = \vec{1} \otimes \ucoeff^n + \Delta t_n (\mat{A} \otimes \mat{I}_\xi) \underline{\vec{F}},
\end{equation}
where the vector $\underline{\ucoeff}^\top = (\ucoeff^1, \dots, \ucoeff^{\Nt})$ contains the $\Nt$ intermediate solution stages and $\underline{\vec{F}}^\top = (\vec{F}(t_n + \Delta t_n c_1, \ucoeff^1), \dots, \vec{F}(t_n + \Delta t_n c_{\Nt}, \ucoeff^{\Nt}))^\top$. Here, $\ucoeff^n$ denotes the RK solution in the previous time step. The new solution is given by
\begin{equation} \label{eq:RK_solution}
\ucoeff^{n+1} = \ucoeff^n + \Delta t_n ( \vec{b}^\top \otimes \mat{I}_\xi ) \underline{\vec{F}}.
\end{equation}

Herein we consider a particular family of implicit RK methods, namely Lobatto IIIC \cite{Jay2015, lehsten:21}. These methods are A-, L- and B-stable and are thus suitable for stiff and nonlinear problems. The order of the $\Nt$-stage Lobatto IIIC method is $2(\Nt-1)$ and the order of the individual stages is $\Nt-1$. 
However, this choice of method is also motivated by its equivalence to a space-time DG-SEM formulation, which is described in the next section. The Butcher tableaus for the 2-, 3- and 4-stage Lobatto IIIC methods are found in \nameref{Appendix_Butcher}.

The following code is an example how to use the Lobatto IIIC solvers in \assimulo: \\

\begin{python}
# import solver form assimulo
import assimulo.ode as aode
import assimulo.solvers as aso
# import Lobatto IIIC solvers
from Lobatto_IIIC_2s import Lobatto2ODE
from Lobatto_IIIC_3s import Lobatto3ODE
from Lobatto_IIIC_4s import Lobatto4ODE

# set up explicit problem, user-defined rhs
prob = aode.Explicit_Problem(rhs, y0, t0) 
# user-defined Jacobian
prob.jac = jacobian 
# choose solver
solver = Lobatto2ODE(prob) 
# run solver until endTime
t, y = solver.simulate(endTime) 
\end{python}


\section{Space-Time DG-SEM} \label{ch:STDG}

We now consider DG-SEM applied to \eqref{eqn:general} with the time variable $t$ treated simply as an additional dimension. The result is space-time DG-SEM.

Defining the gradient $\underline{\nabla} := \left( \nabla, \frac{\partial}{\partial t} \right)$ and the new convective and viscous fluxes
$$
\stFc = \begin{bmatrix}
\Fc & \sol
\end{bmatrix},
\quad
\stFv = \begin{bmatrix}
\Fv & 0
\end{bmatrix},
$$
we can rewrite \eqref{eqn:general} as a $d+1$-dimensional problem over the space-time domain $\underline{\Omega} := \Omega \times (0,T) \subset \RRR^{d+1}$ as
\begin{eqnarray}
\label{eqn:spacetime}
  \underline{\nabla} \cdot\big( \stFc(\sol) -
       \stFv(\sol,\nabla\sol) \big) = S(\sol) \ \ \mbox{ in } \underline{\Omega}.
\end{eqnarray}

Given a tessellation $\underline{\grid}$ of $\underline{\Omega}$ we introduce the piecewise polynomial space
\begin{equation}
\label{eqn:stvspace}
   \underline{\phispace} = \{\underline{\vecv}\in L^2(\underline{\Omega},\RRR^{r}) \; \colon
    \underline{\vecv}|_{\elem}\in[\mathcal{P}_p(\elem)]^{r}, \ \elem\in\underline{\grid}\},
      \ \;p \in \NNN.
\end{equation}
Then the space-time DG-SEM discretization of \eqref{eqn:spacetime} follows analogously to \eqref{eqn:elementint} and \eqref{eqn:surfaceint}:
\begin{equation}
\label{stconvDiscr}
\stdual { \stspcoper(\stdf) } := \stdual{ \underline{K}_h(\stdf) } + \stdual{ \underline{I}_h(\stdf) },
\end{equation}
with the element integrals
%
\begin{eqnarray}
\label{eqn:stelementint}
   \stdual{ \underline{K}_h(\stdf) } &:=&
      \sum_{\elem \in \underline{\grid}} \int_{\elem}
      \big( ( \stFc(\stdf) - \stFv(\stdf, \underline{\nabla} \stdf ) ) : \nabla\stbasefct + \stsu
      \cdot \stbasefct \big) \, dx,
\end{eqnarray}
and the surface integrals
\begin{eqnarray}
\label{eqn:stsurfaceint}
   \stdual{ \underline{I}_h(\df) } &:=&
      \sum_{e \in \underline{\Gamma}_i} \int_e \big(
      \vaver{\stFv(\stdf, \vjump{\stdf} )^\top : \nabla\stbasefct} +
      \vaver{\stFv(\stdf, \underline{\nabla}\stdf)} : \vjump{\stbasefct} \big) \, d\underline{S} \nonumber \\
    &&- \sum_{e \in \underline{\Gamma}} \int_e \big( \stfluxF(\df) - \stfluxA(\stdf,\underline{\nabla}\stdf)\big) :
      \vjump{\stbasefct} \, dS.
\end{eqnarray}
Here, $\underline{\Gamma}_i$ and $\underline{\Gamma}$ have analogous meanings to their spatial counterparts $\Gamma_i$ and $\Gamma$. The numerical fluxes are given by
$$
\stfluxF = \begin{bmatrix}
\fluxF & \sol^*
\end{bmatrix},
\quad 
\stfluxA = \begin{bmatrix}
\fluxA & 0
\end{bmatrix},
$$
where $\sol^*$ is a simple upwind flux in time.

In our considered framework, \eqref{stconvDiscr} can be implemented quite nicely by increasing the
dimension and applying the above discussed modifications\footnote{Note that for the 4D version ($3d + time$) 
a UFL patch (see \nameref{Appendix_UFL}) was added to introduce 
the 4D reference elements to UFL code.}. \\

\begin{python}
d = 2 # 1,2,3 is the spatial dimension
from dune.grid import cartesianDomain, structuredGrid as leafGrid
t_end, timeSteps = 1.0, 10
dt = t_end / timeSteps 
# create grid that tessellates $[0,1]^d \times [0,\Delta t]$ with 10 elements in space and 1 element in time
T_h = leafGrid(cartesianDomain([0]*d + [0], [1]*d + [dt], [10]*d + [1])) # create a space-time grid 
p = 3 # polynomial degree
# create DG space with Lagrange basis and Gauss-Lobatto interpolation points
V_h = dglagrangelobatto( T_h, order=p )

def appendTime( F, u ):  
  return ufl.as_tensor([ *[[F[k,i] if i<d else u[k] for i in range(d+1)] for k in range(len(u))] ])

def F_c( u ):
  from molspacediscr import F_c # import $F_c$ used in MOL discretization
  F_spc = F_c(u) # compute spatial fluxes
  # append time derivative as last column 
  return appendTime( F_spc, u )

def F_v( u ):  
  from molspacediscr import F_v # import $F_v$ used in MOL discretization
  F_spc = F_v(u) # compute spatial fluxes
  # append column of zeros since there is no diffusion in time
  return appendTime( F_spc, [0.]*len(u) )

# trial and test function
u   = TrialFunction(V_h)
psi = TestFunction(V_h)
# element integral from equation \eqref{eqn:stelementint}
K_h = inner(F_c(u) - F_v(u)*grad(u)), grad(psi)) * dx \  # fluxes
    + inner(S(u), psi) * dx                              # source term

# normal and mesh width
n    = FacetNormal(V_h)
h_e  = avg( CellVolume(V_h) ) / FacetArea(V_h)
# penalty parameter for Symmetric Interior Penalty scheme
eta  = Constant( 10*V_h.order**2 if V_h.order > 0 else 1, "penalty" )
# surface integral from equation \eqref{eqn:stsurfaceint}
I_h = inner(jump(H_c(u), jump(psi)) * dS \  # interior skeleton for convective part
    + H_cb(u)*psi*ds \                      # domain boundary for convective part
    - inner(jump(F_v(u),n),avg(grad(psi)))   * dS \   # symmetry term
    - inner(avg(F_v(u)*grad(u)),jump(psi,n)) * dS \   # consistency term
    + eta/h_e*inner(jump(u, avg(F_v(u))*n),jump(psi,n)) * dS  # penalty term
\end{python}

\begin{remark}
It is of practical interest to generalize the space $\underline{\phispace}$ so that the time dimension may be discretized by polynomials of a different order than the spatial dimensions. We will henceforth refer to the number of temporal nodes in each element as $\Nt$ so that the polynomial degree in time is $\Nt-1$. This notation contrasts standard DG terminology, where nodes are typically indexed from $0$ to $p$. Additionally, note that this is the same notation used for the number of stages of the Lobatto IIIC method in Section \ref{ch:MOL}. Stages are typically indexed from $1$ to $s$. However, to minimize the use of notation and to make the connection between the two viewpoints clearer, we write $\Nt$ to count the degrees of freedom within a time element, whether this pertains to the DG or RK interpretation.
\end{remark}

After space-time discretization, the discrete solution $\stdf \in \underline{\phispace}$ takes the form $\stdf(t,x) = \sum_{i,n} \sol_i^n \basefct_i(x) \tbasefct_n(t)$. Here, the sum is taken over all tensor product LGL nodes in $d+1$ dimensions. The vector of coefficients is now given by
\begin{equation}
\label{eq:StCoeff}
\stucoeff = (\ucoeff^1,\dots,\ucoeff^{\Nt})^\top,
\end{equation}
where $\ucoeff^i$ contains all the spatial unknowns in the $i$th time element.

The space-time discretization \eqref{stconvDiscr} can alternatively be derived by starting from \eqref{eqn:ode} and discretizing in time with DG-SEM. Multiplying \eqref{eqn:ode} by a test function $\tbasefct(t)$ and integrating over the $n$th time element results in
$$
\int_{t_n}^{t_{n+1}} \ucoeff_t \tbasefct \text{d}t = \int_{t_n}^{t_{n+1}} \vec{F}(t,\ucoeff(t)) \tbasefct \text{d}t.
$$
We transform this equation to the reference element $[-1,1]$ using the mapping $t=t_n + \frac{\Delta t_n}{2}(1+\tau)$, where $\Delta t_n = t_{n+1} - t_n$. After integration by parts the resulting equation reads
$$
[\ucoeff \tbasefct]_{-1}^{1} - \int_{-1}^{1} \ucoeff \tbasefct_{\tau} \dt = \frac{2}{\Delta t_n} \int_{-1}^{1} \vec{F}(\tau,\ucoeff(\tau)) \tbasefct \dt.
$$
We now follow the steps of DG-SEM, i.e. approximating $\ucoeff$ and $\vec{F}$ by interpolants
\begin{align*}
\ucoeff \approx \sum_{j=1}^{\Nt} \ucoeff^j \tbasefct_j(\tau), \\
\vec{F} \approx \sum_{j=1}^{\Nt} \vec{F}^j \tbasefct_j(\tau),
\end{align*}
and the integrals by Gauss-Lobatto quadrature with nodes $\tau_j$ and weights $\omega_j$. Using the cardinal property of the Lagrange basis polynomials $\tbasefct$, the resulting DG-SEM discretization becomes
\begin{equation} \label{eq:DGSEM_time_component}
\delta_{i \Nt} \ucoeff^{*} - \delta_{i 1} \ucoeff^{*} - \sum_{j=1}^{\Nt} \omega_j \ucoeff^j \left. \frac{\text{d} \tbasefct_i}{\text{d} \tau} \right|_{\tau_j} = \frac{2}{\Delta t_n} \omega_i \vec{F}^i, \quad i=1,\dots,\Nt.
\end{equation}
Here, we have replaced the boundary terms with numerical fluxes $\ucoeff^*$.
With DG-SEM in time, the numerical flux $\underline{\ucoeff}^*$ is always chosen as the upwind flux
\begin{equation} \label{eq:upwind_flux}
\underline{\ucoeff}^* = (\ucoeff^n, \vec{0}, \dots, \vec{0}, \ucoeff^{\Nt})^\top,
\end{equation}
where $\ucoeff^n$ is the numerical solution from the previous time element. This choice leads to an entropy stable numerical scheme if the spatial terms are handled appropriately \cite{Friedrich2019}. It also has the advantage of decoupling the temporal elements. Thus, \eqref{eq:DG_in_time} can be solved as a stand-alone nonlinear system on the $n$th time element.

Defining the boundary, mass and differentiation matrices
\begin{align*}
&\mat{B}_{\tau} = \mathrm{diag}([-1,0,\dots,0,1]) \in \RRR^{N_{\tau} \times N_{\tau}},\\
&\mat{M}_{\tau} = \mathrm{diag}([\omega_1, \dots, \omega_{N_{\tau}}]) \in \RRR^{N_{\tau} \times N_{\tau}},\\
&(\mat{D}_{\tau})_{ji} = \left. \frac{\text{d} \tbasefct_i}{\text{d} \tau} \right|_{\tau_j} \in \RRR^{N_{\tau} \times N_{\tau}},
\end{align*}
we can write \eqref{eq:DGSEM_time_component} in matrix form on each reference element as
\begin{equation} \label{eq:DG_in_time}
\left( \mat{B}_\tau \otimes \mat{I}_\xi \right) \underline{\ucoeff}^*
- \left( \mat{D}_\tau^\top \mat{M}_\tau \otimes \mat{I}_\xi \right) \underline{\ucoeff} = \frac{\Delta t_n}{2} (\mat{M}_\tau \otimes \mat{I}_\xi) \underline{\vec{F}}(\underline{\ucoeff}),
\end{equation}
where $\mat{M}_\tau$ is the local temporal mass matrix and $\mat{M}_\tau \mat{D}_\tau$ defines the corresponding stiffness matrix. Here, $\underline{\vec{F}}^\top(\underline{\ucoeff}) = (\vec{F}^\top(t_n + \frac{\Delta t_n}{2}(1 + \tau_1), \ucoeff^1), \dots, \vec{F}^\top(t_n + \frac{\Delta t_n}{2}(1 + \tau_{N_\tau}), \ucoeff^{N_\tau}))$, where $\underline{\ucoeff}$ is given by \eqref{eq:StCoeff} and $\tau_k$ is the $k$th LGL node; see \cite{Gassner2013} for details. The operation $\otimes$ denotes the Kronecker product and $\mat{I}_\xi$ is the identity matrix whose dimension is given by the number of spatial nodes. 

We finish this section by remarking that while \eqref{stconvDiscr} describes the global space-time DG-SEM discretization, the alternative formulation \eqref{eq:DG_in_time} pertains to a single time element.


\section{Theoretical Aspects of Space-Time DG-SEM} \label{ch:theory}
In this section we discuss important properties of the space-time DG-SEM, in particular the equivalence of the temporal discretization and the Lobatto IIIC family of Runge-Kutta methods. To make the connection between DG-SEM and Runge-Kutta methods clear, we consider the solution at the final point in the time element, i.e.
\begin{equation} \label{eq:un}
\ucoeff^{\Nt} \equiv (\vec{e}_{\Nt}^\top \otimes \mat{I}_\xi) \underline{\ucoeff},
\end{equation}
where $\vec{e}_{\Nt}^\top = (0,\dots,0,1) \in \RRR^{\Nt}$. We will also make use of the vector $\vec{e}_1^\top = (1,0,\dots,0) \in \RRR^{\Nt}$. Following \cite{Boom2015}, we set out to show that $\ucoeff^{\Nt} = \ucoeff^{n+1}$, where $\ucoeff^{n+1}$ is the numerical solution arising from the Lobatto IIIC method in \eqref{eq:RK_solution}, under the assumption that this equality holds in the previous (i.e. in the $(n-1)$st) time element.  

The DG-SEM discretization \eqref{eq:DG_in_time} constitutes a so called \emph{Summation-By-Parts} (SBP) method \cite{Gassner2013}, meaning that the following conditions are satisfied:
\begin{equation} \label{eq:SBP}
\mat{M}_{\tau} = \mat{M}_{\tau}^\top > \mat{0}, \quad \mat{M}_{\tau} \mat{D}_{\tau} + (\mat{M}_{\tau} \mat{D}_{\tau})^\top = \mat{B}_{\tau}.
\end{equation}
The SBP property \eqref{eq:SBP} is at the heart of the connection of DG-SEM in time to implicit Runge-Kutta methods.

\subsection{DG-SEM and Lobatto IIIC}

SBP methods were historically developed to be used as spatial discretizations \cite{Kreiss1974,Strand1994}. For an overview of these techniques, see \cite{Svard2014,Fernandez2014}. In recent years, their use as time stepping schemes has been explored \cite{Nordstrom2013} and connections to implicit Runge-Kutta methods have been discovered \cite{Boom2015}. Here we summarize the steps showing that \eqref{eq:DG_in_time} can be reformulated as an implicit RK method applied to the system of ODEs \eqref{eqn:ode}.

We begin by using the SBP property \eqref{eq:SBP} in the second term of \eqref{eq:DG_in_time} and then multiplying by $(\mat{M}_\tau^{-1} \otimes \mat{I}_\xi)$ to obtain the so called \emph{strong form},
\begin{equation} \label{eq:DG_in_time_strong_form}
(\mat{D}_\tau \otimes \mat{I}_\xi) \underline{\ucoeff} = (\mat{M}_\tau^{-1} \mat{B}_\tau \otimes \mat{I}_\xi) (\underline{\ucoeff} - \underline{\ucoeff}^*) + \frac{\Delta t_n}{2} \underline{\vec{F}}.
\end{equation}
Note that $\mat{B}_\tau = \vec{e}_{N_\tau} \vec{e}_{N_\tau}^\top - \vec{e}_1 \vec{e}_1^\top$ and $(\vec{e}_{N_\tau}^\top \otimes \mat{I}_\xi )(\underline{\ucoeff} - \underline{\ucoeff}^*) = (\ucoeff^{\Nt} - \ucoeff^{\Nt}) = \vec{0}$. Using \eqref{eq:upwind_flux}, the second term in \eqref{eq:DG_in_time_strong_form} can therefore be expressed as
$$
(\mat{M}_\tau^{-1} \mat{B}_\tau \otimes \mat{I}_\xi) (\underline{\ucoeff} - \underline{\ucoeff}^*) = -(\mat{M}_\tau^{-1} \otimes \mat{I}_\xi) [ (\vec{e}_1 \vec{e}_1^\top \otimes \mat{I}_\xi) \underline{\ucoeff} - (\vec{e}_1 \otimes \ucoeff^n)].
$$
Grouping together terms that multiply the solution $\underline{\ucoeff}$, we rewrite \eqref{eq:DG_in_time_strong_form} as
\begin{equation} \label{eq:DG_D_tilde}
((\mat{D}_\tau + \mat{M}_\tau^{-1} \vec{e}_1 \vec{e}_1^\top) \otimes \mat{I}_\xi) \underline{\ucoeff} = (\mat{M}_\tau^{-1} \vec{e}_1 \otimes \ucoeff^n) + \frac{\Delta t_n}{2} \underline{\vec{F}}.
\end{equation}
Next, we multiply \eqref{eq:DG_D_tilde} by $((\mat{D}_\tau + \mat{M}_\tau^{-1} \vec{e}_1 \vec{e}_1^\top) \otimes \mat{I}_\xi)^{-1}$. Upon doing this, first note that
$$
(\mat{D}_\tau + \mat{M}_\tau^{-1} \vec{e}_1 \vec{e}_1^\top)^{-1} \mat{M}_\tau^{-1} \vec{e}_1 = \vec{1} := (1,\dots,1)^\top \in \RRR^{\Nt},
$$
which follows from observing that $(\mat{D}_\tau + \mat{M}_\tau^{-1} \vec{e}_1 \vec{e}_1^\top) \vec{1} = \mat{M}_\tau^{-1} \vec{e}_1$ since $\mat{D}_\tau \vec{1} = \vec{0}$ by consistency. Thus, the following system arises:
\begin{equation} \label{eq:DG_in_RK_form}
\begin{aligned}
\underline{\ucoeff} &= \vec{1} \otimes \ucoeff^n + \Delta t_n \frac{1}{2} ((\mat{D}_\tau + \mat{M}_\tau^{-1} \vec{e}_1 \vec{e}_1^\top) \otimes \mat{I}_\xi)^{-1} \underline{\vec{F}} \\
&= \vec{1} \otimes \ucoeff^n + \Delta t_n \frac{1}{2} ((\mat{D}_\tau + \mat{M}_\tau^{-1} \vec{e}_1 \vec{e}_1^\top)^{-1} \otimes \mat{I}_\xi) \underline{\vec{F}}.
\end{aligned}
\end{equation}

The equation system \eqref{eq:DG_in_RK_form} should be compared with the stage equations \eqref{eq:RK_stage_equations} that arose from the MOL discretization using implicit RK. We see that the temporal DG-SEM discretization defines an RK method with coefficient matrix $\mat{A} = \frac{1}{2} (\mat{D}_\tau + \mat{M}_\tau^{-1} \vec{e}_1 \vec{e}_1^\top)^{-1}$ and nodes $\vec{c} = (\vec{1} + \vec{\tau})/2$, where $\vec{\tau} = (\tau_1, \dots, \tau_{\Nt})^\top$ is the vector of LGL nodes. Further, the vector $\underline{\ucoeff}$, which in the DG-SEM context contains the interpolation coefficients $u_i^n$, has adopted the role of the stage vectors of the RK method.

To complete the transition from DG-SEM to RK, we compute the numerical solution at the final time node, $\ucoeff^{\Nt} = (\vec{e}_{N_\tau}^\top \otimes \mat{I}_\xi) \underline{\ucoeff}$. To this end we observe that the SBP property \eqref{eq:SBP} gives the relation
$$
\vec{1}^\top \mat{M}_\tau (\mat{D}_\tau + \mat{M}_\tau^{-1} \vec{e}_1 \vec{e}_1^\top) = \vec{1}^\top (\vec{e}_{N_\tau} \vec{e}_{N_\tau}^\top - \mat{D}_\tau^\top \mat{M}_\tau) = \vec{e}_{N_\tau}^\top,
$$
so that
$$
\vec{e}_{N_\tau}^\top (\mat{D}_\tau + \mat{M}_\tau^{-1} \vec{e}_1 \vec{e}_1^\top)^{-1} = \vec{1}^\top \mat{M}_\tau.
$$
Consequently, multiplying \eqref{eq:DG_in_RK_form} by $(\vec{e}_{N_\tau}^\top \otimes \mat{I}_\xi) \underline{\ucoeff}$ yields
\begin{equation} \label{eq:DG_RK_solution}
\ucoeff^{\Nt} = \ucoeff^n + \Delta t_n \frac{1}{2} ( \vec{1}^\top \mat{M}_\tau \otimes \mat{I}_\xi ) \underline{\vec{F}}.
\end{equation}
Comparing \eqref{eq:DG_RK_solution} with the solution \eqref{eq:RK_solution} of the implicit RK method, we see that the vector $\vec{b}$ in the Butcher tableau is related to the DG-SEM discretization by $\vec{b}^\top = \vec{1}^\top \mat{M}_\tau/2$, and that the RK solution is simply the $N_\tau$th component of the DG solution $\underline{\ucoeff}$.

To summarize, a DG-SEM time discretization is equivalent to an implicit RK method whose Butcher tableau is defined in terms of the DG method as
\begin{equation} \label{eq:SBP_to_RK}
\mat{A} = \frac{1}{2}(\mat{D}_\tau + \mat{M}_\tau^{-1} \vec{e}_1 \vec{e}_1^\top)^{-1}, \quad \vec{b} = \frac{1}{2} \mat{M}_\tau \vec{1}, \quad \vec{c} = \frac{\vec{1} + \vec{\tau}}{2}.
\end{equation}
The two methods yield two different nonlinear systems; for DG-SEM and RK they are respectively given by
\begin{subequations}
\begin{align}
\left( \mat{B}_\tau \otimes \mat{I}_\xi \right) \underline{\ucoeff}^*
- \left( \mat{D}_\tau^\top \mat{M}_\tau \otimes \mat{I}_\xi \right) \underline{\ucoeff} &= \frac{\Delta t_n}{2} (\mat{M}_\tau \otimes \mat{I}_\xi) \underline{\vec{F}}, \label{eq:DG_vs_RK:DG} \\
\underline{\ucoeff} &= \vec{1} \otimes \ucoeff^n + \Delta t_n (\mat{A} \otimes \mat{I}_\xi) \underline{\vec{F}}. \label{eq:DG_vs_RK:RK}
\end{align}
\end{subequations}
These systems have the same solution $\underline{\ucoeff}$ since we can transition from \eqref{eq:DG_vs_RK:DG} to \eqref{eq:DG_vs_RK:RK} in a series of algebraic steps. More precisely, the connection is made by rewriting \eqref{eq:DG_vs_RK:DG} in strong form, then multiplying by $(\mat{D}_\tau + \mat{M}_\tau^{-1} \vec{e}_1 \vec{e}_1^\top \otimes \mat{I}_\xi)^{-1}$.

Note that the latter step demands that $\mat{D}_\tau + \mat{M}_\tau^{-1} \vec{e}_1 \vec{e}_1^\top$ is invertible. This is the case if and only if $\mat{D}_\tau$ is \emph{null-space consistent}, i.e. if $\mathrm{ker}(\mat{D}_\tau) = \mathrm{span}(\vec{1})$ \cite{Linders2020}. This is known to hold for all $N_\tau > 1$ \cite{Ruggiu2018,Linders2022}.

Finally, the Butcher tableau formed from \eqref{eq:SBP_to_RK} coincides with that of the Lobatto IIIC family of implicit Runge-Kutta methods. This follows from the use of LGL nodes and quadrature weights, together with a set of accuracy conditions satisfied by the two formulations \cite{Ranocha2019}. We will detail these in the next section. The derivation above therefore shows that DG-SEM in time and the Lobatto IIIC methods are mathematically equivalent, and that we in fact have $\ucoeff^{\Nt} = \ucoeff^{n+1}$. The coefficients for the DG-SEM matrices with $\Nt \in \{2,3,4\}$ are listed in \nameref{Appendix_DGSEM}.

\subsection{Comparison of terminology}

While DG-SEM in time and Lobatto IIIC are algebraically equivalent methods, they have been developed in different research communities and disparities have consequently arisen in terms of terminology. This pertains in particular to the notions of order and stability.

Beginning with RK methods, we take as our starting point the system of ODEs \eqref{eqn:ode}. The (classical) notion of order is defined as follows:

\begin{definition} \label{def:RK_order}
A Runge-Kutta method is of \emph{order} $p$ if
$$
\| \ucoeff^{n+1} - \ucoeff(t_{n+1}) \| =  \mathcal{O}(\Delta t^{p}), \quad \Delta t \rightarrow 0
$$
holds, whenever problem \eqref{eqn:ode} is sufficiently smooth.
\end{definition}
The norm can be any vector norm and $\ucoeff^{n+1} - \ucoeff(t_{n+1})$ is called the global error. 

The classical order of RK methods is determined by certain order conditions. To make the connection with DG-SEM as clear as possible, we present here a set of simplified conditions that are sufficient for the method to be of order $p$ \cite{Butcher1964}:

\begin{theorem} \label{thm:order_conditions}
Suppose that an implicit Runge-Kutta method satisfies the conditions
\begin{description}
\item[$B(p_B)$:] $\vec{b}^\top \vec{c}^{j-1} = \frac{1}{j}, \quad j = 1, \dots, p_B$,
\item[$C(p_C)$:] $\mat{A} \vec{c}^{j-1} = \frac{\vec{c}^j}{j}, \quad j = 1, \dots, p_C$,
\item[$D(p_D)$:] $\mat{A}^\top \text{diag}(\vec{b}) \vec{c}^{j-1} = \frac{1}{j} \text{diag}(\vec{b}) (\vec{1} - \vec{c}^j), \quad j = 1, \dots, p_D$,
\end{description}
where $p_B \leq 2(p_C+1)$ and $p_B \leq p_C + p_D + 1$. Then the method is of order $p = p_B$.
\end{theorem}

The conditions $C(p_C)$ play a particularly important role in the context of stiff problems and have its own moniker:

\begin{definition}
A Runge-Kutta method that satisfies the order conditions $C(p_C)$ is said to have \emph{stage order} $p_C$.
\end{definition}

The stage order of the RK method describes the accuracy with which the intermediate stages are approximated. We will delve into the meanings of the conditions $B$, $C$ and $D$ shortly. However, first we summarize the various order concepts for Lobatto IIIC; see \cite[Chapter IV.5]{Hairer2010}.

\begin{theorem} \label{thm:Lobatto_IIIC_order}
The Lobatto IIIC method with $\Nt$ stages satisfies $B(2\Nt-2)$, $C(\Nt-1)$ and $D(\Nt-1)$. Consequently it has stage order $\Nt-1$ and is of order $2\Nt-2$.
\end{theorem}

We momentarily leave the RK viewpoint and focus on DG methods. DG-SEM was developed for spatial discretizations of time-dependent PDEs. Errors are measured in an $L_2$ norm over a spatial domain. For a discrete solution $\vec{u}_h$, this norm can be computed via the quadrature rule exactly: 
\begin{equation} \label{eq:DG_norm}
\| \vec{u}_h \|_{L_2(\Omega)} =  \sum_i \left( \vec{u}_{h_i}^\top \mat{M}_{\xi_i} \vec{u}_{h_i} \right)^{\frac{1}{2}}.
\end{equation}
The sum is taken over all elements and $\mat{M}_{\xi_i}$ is the local spatial mass matrix on element $i$. Assuming vanishingly small errors from the time discretization, the order of convergence measured in this norm is typically $p+1$ or $p+\frac{1}{2}$, depending on the nature of problem \eqref{eqn:general}, the choice of numerical fluxes, and sometimes on whether $p$ is odd or even \cite{Hesthaven2008,Zhang2006}.

Conversely, when using DG-SEM as a time integration method, one works in the space $L_2([0,T])$. With an upwind flux in time, for sufficiently smooth and nonstiff problems, the order of convergence in this norm is $\Nt$ \cite{Lundquist2014}. This order is much smaller than the one of the Lobatto IIIC method, which requires some discussion. 

There are several other order concepts in the DG literature. Here we follow \cite{Boom2015} and relate these to the corresponding concepts in the RK framework.

\begin{itemize}
\item The \emph{order of the operator} is the highest degree $q$ for which $\mat{D}_\tau \vec{\tau}^q = q \vec{\tau}^{q-1}$. The exponentiation should be interpreted elementwise, and we take $\vec{\tau}^0 = \vec{1}$ as a definition. For DG-SEM we have $q = \Nt-1$.

Multiplying $C(p_C)$ by $\vec{A}^{-1}$ as given in \eqref{eq:SBP_to_RK} and utilizing the fact that the first element in $\vec{c}$ is zero, we see that the RK order condition $C(p_C)$ actually describes precisely the order of the operator $\mat{D}_\tau$. A transformation of the reference element to $[0,1]$ is necessary in this step. In other words, the order of the operator is a concept identical to the stage order of the corresponding Lobatto IIIC method.
\item The \emph{order of the norm/quadrature/mass matrix} is the highest degree $m$ such that $(m+1) \vec{1}^\top \mat{M}_\tau \vec{\tau}^m = 1 - (-1)^{m+1}$, i.e. for which $\mat{M}_\tau$ exactly integrates polynomials. For DG-SEM, $\vec{1}^\top \mat{M}_\tau$ is a row vector with the $\Nt$ weights of the Gauss-Lobatto quadrature rule and we consequently have $m = 2(\Nt-1)$.

Using \eqref{eq:SBP_to_RK} we note that the condition $B(p_B)$ simply describes the order of the quadrature, although applied to $\vec{c}$ rather than $\vec{\tau}$. Again, this amounts to a transformation from $\tau \in [-1,1]$ to $[0,1]$.
\item Pertinently, it turns out that the order of accuracy of the final component $\ucoeff^{\Nt} \equiv \ucoeff^{n+1}$ is $2(\Nt - 1)$ \cite{Lundquist2014}, at least for smooth nonstiff problems. This \emph{superconvergence} can be proven using the theory of dual consistent SBP methods \cite{Hicken2011}. Here it suffices to say that it is a consequence of the order of the quadrature and choosing the upwind numerical flux \eqref{eq:upwind_flux}.

The superconvergence result pertaining to DG-SEM corresponds to the classical order of Lobatto IIIC as introduced in Definition \ref{def:RK_order}. Note that this is a consequence of considering the pointwise error in time rather than $\|\cdot \|_{L_2([0,T])}$.
\end{itemize}

To the best of our knowledge, conditions $D(p_D)$ have no clear interpretation in the language of DG. Nevertheless, using the SBP property \eqref{eq:SBP} and the diagonality of $\mat{M}_\tau$ it is shown in \cite{Boom2015} that $C(p_C)$ is satisfied with $p_C = \Nt-1$, which is consistent with Theorem \ref{thm:Lobatto_IIIC_order}.


The convergence theory for RK methods relies on certain regularity properties of the problem being solved. In particular, it is assumed that the right-hand side of the system of ODEs \eqref{eqn:ode} satisfies a one-sided Lipschitz condition,
\begin{equation} \label{eq:one-sided_Lipschitz}
\langle \ucoeff - \vec{v}, \vec{F}(t,\ucoeff) - \vec{F}(t,\vec{v}) \rangle \leq \beta \| \ucoeff - \vec{v} \|^2,
\end{equation}
where $\langle \cdot, \cdot \rangle$ denotes some inner product and $\| \cdot \|$ the corresponding norm. If $\beta \leq 0$, the problem is \emph{contractive}. DG-SEM in time, and hence Lobatto IIIC, is stable for contractive problems, i.e. they are B-stable methods. Convergence for contractive problems is correspondingly known as B-convergence. B-convergence can be shown for Lobatto IIIC if $\beta < 0$, but in general not if $\beta = 0$ if $\Nt > 2$; see \cite{Schneid1987} for details.

Regularity of the type \eqref{eq:one-sided_Lipschitz} is not standard in the literature on spatial discretization using high order DG methods. Rather, estimates of the form $\langle \ucoeff, \vec{F}(t,\ucoeff) \rangle \leq 0$ are common. Such discretizations are referred to as \emph{semi-bounded}, or in Runge-Kutta parlance, as \emph{monotonic}. They arise particularly for discretizations of linear, homogeneous hyperbolic or parabolic problems, but also e.g. for the velocity components of the incompressible Navier-Stokes equations \cite{Nordstrom2019}.

If $\vec{F}(t,\vec{0}) = \vec{0}$, then semi-boundedness is a special case of contractivity and results on B-stability and B-convergence apply. However, for e.g. the equations of compressible flow, semi-boundedness must typically be replaced by \emph{entropy stability}, i.e. regularity of the form $\langle \eta'(\ucoeff), \vec{F}(t,\ucoeff) \rangle \leq 0$. Here, $\eta$ is some convex function of $\ucoeff$ referred to as an entropy \cite{Fisher2013,Fisher20132}. 
A convergence theory for implicit RK methods applied to entropy stable problems is desirable but currently appears to be missing from the literature.


\section{Practical Aspects of Space-Time DG-SEM} \label{ch:practical}
In this section we discuss two archetypal implementations: On the one hand, the method
of lines approach with Lobatto IIIC as discussed in Section \ref{ch:MOL} and on the other 
hand the space-time DG approach, as discussed in Section \ref{ch:STDG}. We will refer to 
these as LoDG and STDG, respectively.

It is of course possible to produce a code that uses elements from both the LoDG and STDG formulation and 
thereby falls somewhere in between these approaches. 
However, here we adopt the point of view of a user who seeks to 
use an already available code base rather than producing a brand new solver.

Even though STDG and LoDG are mathematically equivalent methods, their respective implementations differ in several key aspects, each with particular requirements and accompanying merits. Here we will outline several such differences, and the choices a user will inevitable face when deciding on which implementation to select.

Several multi-dimensional DG-SEM solvers exist, such as 
Nektar++ \cite{Karniadakis2013}, Fluxo, Flexi and the latest iteration Trixi \cite{ranocha2021adaptive} and others.
In particular, this approach is popular for weather and climate prediction  and has been used e.g. in NUMA \cite{numa} and HOMAM \cite{homam:16}. Thus, in the following discussion we assume that the user has
access to a multi-dimensional DG solver for spatial discretization.

Our work here is based on the \dunefem framework, hence the challenges outlined below are flavoured by this choice. Depending on the software framework at hand,
a user may find that one approach is easier to implement than the other. 
 
\subsection{STDG}

As described in Section \ref{ch:STDG}, the defining feature of STDG is the treatment of the $d$-dimensional time-dependent problem as a $d+1$-dimensional stationary problem. \\

\noindent \textbf{Requirements:} The problem description in the code must be 
extended to a $d+1$-dimensional stationary PDE, which requires the software to be able to handle such problems. In particular, this includes the assembly of mass and stiffness matrices as well as having access to appropriate solvers for the resulting nonlinear system.

This new stationary PDE requires the use of different numerical fluxes 
in space and time; the temporal direction follows a causality principle, 
enforced by the upwind flux, which may not be the best choice for the spatial directions. 

It is desirable to be able to choose different orders of accuracy in space and time, 
which then needs to be made possible in the DG code. In \dune this is implemented for
certain DG spaces \cite{gersbacher} but not yet available for the Lagrange basis
used in this work. 

The numerical solution should be accessible at specific time points in order for the user to visualize
intermediate results and the final solution. 
For $d=3$, this includes extracting $3$-dimensional slices from $4$-dimensional data sets. Related to this issue is the problem of $4$-dimensional mesh generation. An example of how to handle this for the STDG ansatz
in \dunefem is found in Section \ref{ch:STDG}. \\

\noindent \textbf{Merits:} If the requirements above are met, then existing software can be reused to solve the problem, which implies full control over the code. 
Moreover, only one code is needed. This code closely follows the mathematical derivation of the
space-time method and may therefore be more intuitive than alternatives.
 
Due to the relatively simple adaption of an existing code for spatial problems of dimension $d<3$, 
the STDG approach is fast for preliminary testing. An existing DG code is most likely 
optimized for computational resources and might even allow for 
parallelization in time by solving for several time steps at once. 

In summary, this technique allows re-usability and full control over the code.

\subsection{LoDG}

The defining feature of LoDG is the method of lines approach outlined in Section \ref{ch:MOL}. In this DG-SEM solver, each time step is solved individually. \\

\noindent \textbf{Requirements:} 
The (spatial) DG-SEM code needs to be coupled with an ODE solver with an implementation
of a Lobatto IIIC method, most likely coming from another code. Difficulties may arise from the particular requirements of the two codes, such as interfaces for time and space adaptivity, parallelization etc.

The ODE solver might require input in a specific 
format not native to the DG code. Further, an efficient solution procedure may require information from the DG-SEM solver not directly available, such as the Jacobian of the spatial discretization.
 
An example of a Lobatto IIIC solver implemented in \assimulo \cite{Andersson2015} is found at the end of Section \ref{ch:MOL}. \\


\noindent \textbf{Merits:} No adaption of the DG code with respect to the PDE or its dimension is necessary. 
There are no additional difficulties arising in the treatment of $3$-dimensional problems, and the resulting solution can be visualized in a straightforward way. 

Most existing ODE solvers are optimized and equipped with several options, 
for instance adaptive time stepping. Intermediate results are easy to access and 
the order of accuracy in space and time can be chosen independently.

In summary, this technique provides flexibility and allows reuse of existing
simulation workflows.


\subsection{Algorithmic Aspects}

In the following we suppose that we have overcome the most important challenges of the two approaches presented in the previous subsections. Thus we now have access to
\begin{enumerate}[label=(\alph*)]
\item a code that generates a $d+1$-dimensional space-time DG-SEM discretization by following the steps in Section \ref{ch:STDG} (STDG).
\item a code that generates a $d$-dimensional spatial DG-SEM discretization by 
  following the steps in Section \ref{ch:MOL}, and a code for time marching using Lobatto IIIC (LoDG),
\end{enumerate}

In each time step, the LoDG code (approximately) solves \eqref{eq:DG_vs_RK:RK} while the STDG code solves \eqref{eq:DG_vs_RK:DG}, or equivalently solves $\stspcoper(\stdf) = 0$ from \eqref{stconvDiscr}. We now discuss the impact of choosing LoDG or STDG on a variety of algorithmic aspects. \\


\noindent \textbf{Accessing time steps and stages:} Accessing the numerical solution at a particular time is straightforward in most ODE solvers. The solution is computed either by aligning the step sizes with the target times or through accurate interpolation.

While it is possible to implement such techniques with STDG, they are unlikely to be available in a pre-existing DG code. Further, the code will return the numerical solution at all points in one (or several) time elements simultaneously. In fact, it is not obvious that the STDG code will be able to return $\ucoeff^{n+1}$ in a simple way since this requires the extraction of a specific subset of coefficients from the numerical solution vector $\stucoeff$. Yet, this may be necessary e.g. for visualization, to use adaptive time stepping, or in case the solution needs to be filtered or otherwise modified between time steps. The solution can in principle be constructed using $\ucoeff^{n+1} = (\vec{e}_{N_\tau}^\top \otimes \mat{I}_\xi) \underline{\ucoeff}$ with $\vec{e}_{\Nt}^\top = (0,\dots,0,1)$, as was done in Section \ref{ch:theory}. However, this assumes that the ordering of the unknowns in $\stucoeff$ is identical to the one used in that analysis. If not, $(\vec{e}_{N_\tau}^\top \otimes \mat{I}_\xi)$ must be suitably permuted into some matrix $\mat{E}_{\Nt}$ before application. Finding the appropriate permutation matrix may be a nontrivial task, in particular in 4D.

On the other hand, with STDG we have access to all intermediate time stages by default, something that may be challenging with LoDG. This may be useful to compute $L_2$ errors of the numerical solution and has the additional advantage of allowing visualization of the solution away from the main time steps. \\

\noindent \textbf{Adaptive time-stepping:} Adaptive time-stepping for RK methods is standard in modern software and thus will be available in an implementation of LoDG. It requires a way of estimating the numerical error in the next time step. This information is used to adapt the time step to fit a predefined tolerance. Embedding techniques use a vector $\hat{\vec{b}}$ to compute a second numerical solution $\hat{\ucoeff}^{n+1}$ from \eqref{eq:RK_solution} whose accuracy is one order lower than that of $\ucoeff^{n+1}$. The difference $\ucoeff^{n+1} - \hat{\ucoeff}^{n+1}$ can be used to estimate the local error without the need to solve the nonlinear system \eqref{eq:DG_vs_RK:RK} more than once. A detailed strategy for estimating the error and choosing the time step based on the embedding technique is available in \cite[Chapter IV.8]{Hairer2010} for the Radau IIA method, but can be easily adapted to Lobatto IIIC \cite{lehsten:21}. .

An STDG code that follows the steps outlined in Section \ref{ch:STDG} will not have an embedded method. However, if the matrix $\mat{E}_{\Nt}$ can be found that extracts $\ucoeff^{n+1}$, then it is also possible to construct a matrix $\hat{\mat{E}}_{\Nt}$ that extracts $\hat{\ucoeff}^{n+1}$ such that an embedding technique can be used. This procedure may be more invasive than desirable.

Alternatively, the numerical error can be estimated using Richardson extrapolation \cite[Chapter II.4]{Hairer2009}. This procedure requires solving the nonlinear system \eqref{eq:DG_vs_RK:DG} three times; once with a step size $2\Delta t$ and twice with a step size $\Delta t$. The difference between the two solutions yields an error estimate. Due to its expense, this approach hardly seems feasible for a 4D problem. \\

\noindent \textbf{Adaptive Mesh Refinement (AMR):} Since STDG uses a 4D mesh it is straightforward to set up a system that accounts for multiple temporal elements at once, which is not possible with LoDG. This introduces the possibility of using AMR in time in addition to space; see \cite{Jayasinghe2018,Chen2008} and the references therein. Like with Richardson extrapolation, using AMR forces us to solve the nonlinear system multiple times. Additionally, the system now consists of multiple coupled time steps. However, the additional cost may be offset by two factors: Firstly, we expect that the number of degrees of freedom necessary to achieve a given accuracy is significantly reduced by the AMR. Secondly, parallelism can be employed in the temporal direction.

Space-time AMR is not likely to be simple to set up with commercially available software. However, if the initial hurdles can be circumvented, then it is in principle possible to use completely unstructured space-time grids with h/p-refinement. The technique requires a generator for unstructured cuboid meshes in 4D (tesseracts) \cite{Caplan2020} and a way of estimating the numerical error in the final time. Such tools have been developed for 4D simplex meshes in \cite{Yano2012,Caplan2020}, but appear to be missing for other mesh types. \\

\noindent \textbf{Shock capturing and limiting:}  
The DG spatial discretizations used with RK time stepping are 
stable when applied to linear problems such as linear hyperbolic systems.
However, for nonlinear problems spurious oscillations occur near strong shocks or steep gradients.
In this case the DG method requires some extra stabilization unless a first
order scheme ($p=0$) is used that produces a monotonic structure in the shock region.
For higher order schemes many approaches have been suggested 
to make this property available without introducing an excessive amount of
numerical viscosity, which is a characteristic feature of first order schemes.
Several approaches exist, including slope 
limiters, artificial diffusion (viscosity) techniques, and even a posteriori techniques and order reduction methods. A comprehensive literature list is presented in \cite{Shu-boundpreserving-review:16}. 

In \dunefemdg \cite{dunefemdg:21}, both limiter based approaches and artificial diffusion are available to stabilize a 
DG scheme. The slope limiter based approach implemented in \dunefemdg is coupled with a troubled cell
indicator which makes the overall scheme highly non-linear and therefore
not suitable for implicit methods, since the selection of troubled cells could
change between linear iterations and lead to divergence of the linear solver. On
the other hand, artificial diffusion approaches require a discretization of a diffusion term.
Stabilization diffusion coefficients only need to 
be re-computed every time step. A standard approach is available in \dunefemdg. 
Thus, it is more suitable to apply artificial diffusion techniques 
for problems where strong shocks occur. \\

\noindent \textbf{Nonlinear solvers and preconditioning:} The solutions of the nonlinear systems \eqref{eq:DG_vs_RK:DG} and \eqref{eq:DG_vs_RK:RK} must be approximated, typically by iterative methods for large systems. There have been a multitude of suggestions on how to design such methods; early solvers for implicit Runge-Kutta methods based on modified Newton iterations were introduced in \cite{Butcher1976} and \cite{Bickart1977}. A more optimized algorithm is described in \cite{Hairer2010}, and many later developments use this as a starting point. These can be considered black box solvers in the sense that they do not utilize information about the spatial terms in the solution process.

Methods designed specifically for spatial DG discretizations and implicit Runge-Kutta methods are found in e.g. \cite{Pazner2017,Pazner2018I}. Likewise, nonlinear solvers designed for space-time DG and FEM discretizations have been developed \cite{Klaij2007,Sudirham2006}.

Unless the user is willing to make the (possibly considerable) effort to develop and/or implement a nonlinear solver specifically designed for LoDG or STDG, the natural recourse is to use a black box solver. Efficiency gains can possibly be made by introducing a preconditioner designed for DG discretizations; see e.g. \cite{Versbach2019,Pazner2018I,Kasimir2021} for recent developments. However, attention must be payed to the fact that the systems \eqref{eq:DG_vs_RK:DG} and \eqref{eq:DG_vs_RK:RK} have different algebraic properties and therefore likely will respond differently to preconditioners and nonlinear solvers.

For the nonlinear LoDG system \eqref{eq:DG_vs_RK:RK}, the Jacobian is given by
\begin{align} \label{eq:Jacobian_A}
\underline{\mat{I}} - \Delta t_n (\mat{A} \otimes \mat{I}_\xi) \mathcal{J}(\underline{\vec{F}}),
\end{align}
where $\mathcal{J}(\underline{\vec{F}})$ contains the Jacobian of the spatial discretization. The solver in \cite{Hairer2010}, and many recent developments that build upon it, instead use the mathematically equivalent
\begin{align} \label{eq:Jacobian_Ainv}
(\Delta t_n \mat{A})^{-1} \otimes \mat{I}_\xi - \mathcal{J}(\underline{\vec{F}}).
\end{align}
For the nonlinear STDG system \eqref{eq:DG_vs_RK:DG}, the Jacobian is given by
\begin{align} \label{eq:Jacobian_ST}
\left( \mat{D}_\tau^\top \mat{M}_\tau - \vec{e}_{\Nt} \vec{e}_{\Nt}^\top \right) \otimes \mat{I}_\xi + \frac{\Delta t_n}{2} (\mat{M}_\tau \otimes \mat{I}_\xi) \mathcal{J}(\underline{\vec{F}}).
\end{align}
Note from \eqref{eq:SBP_to_RK} that \eqref{eq:Jacobian_Ainv} arises by multiplying \eqref{eq:Jacobian_A} by $(\mat{D}_\tau + \mat{M}_\tau^{-1} \vec{e}_1 \vec{e}_1^\top) \otimes \mat{I}_\xi$. This formulation is therefore very closely related to the STDG Jacobian \eqref{eq:Jacobian_ST}. In fact, they only differ by an application of the SBP property \eqref{eq:SBP} and a multiplication by the temporal mass matrix.

\begin{figure}[h!]
\centering
\begin{subfigure}[b]{.475\textwidth}
\centering
\includegraphics[width=\textwidth]{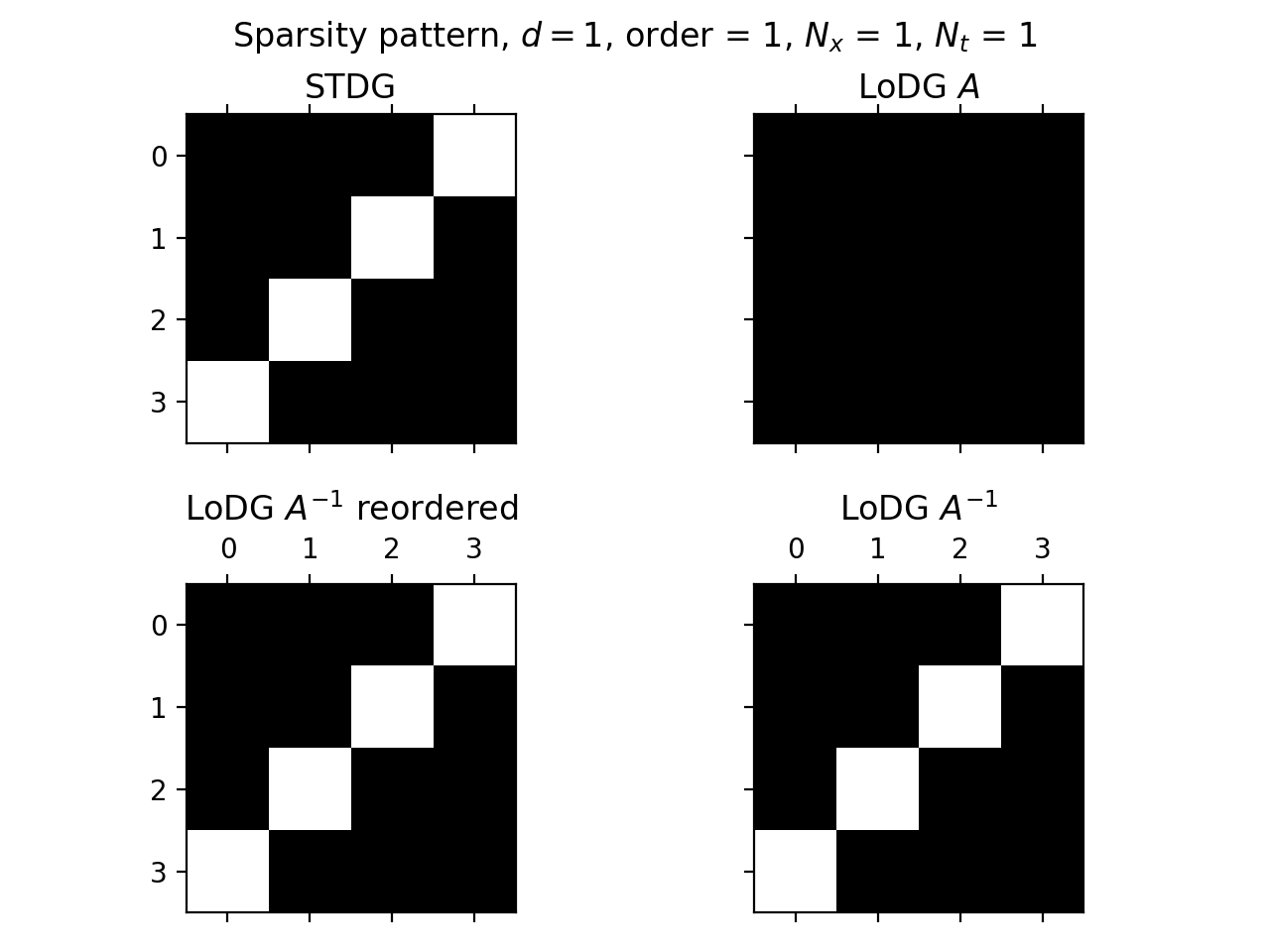}
\end{subfigure}
%
%
\begin{subfigure}[b]{.475\textwidth}
\centering
\includegraphics[width=\textwidth]{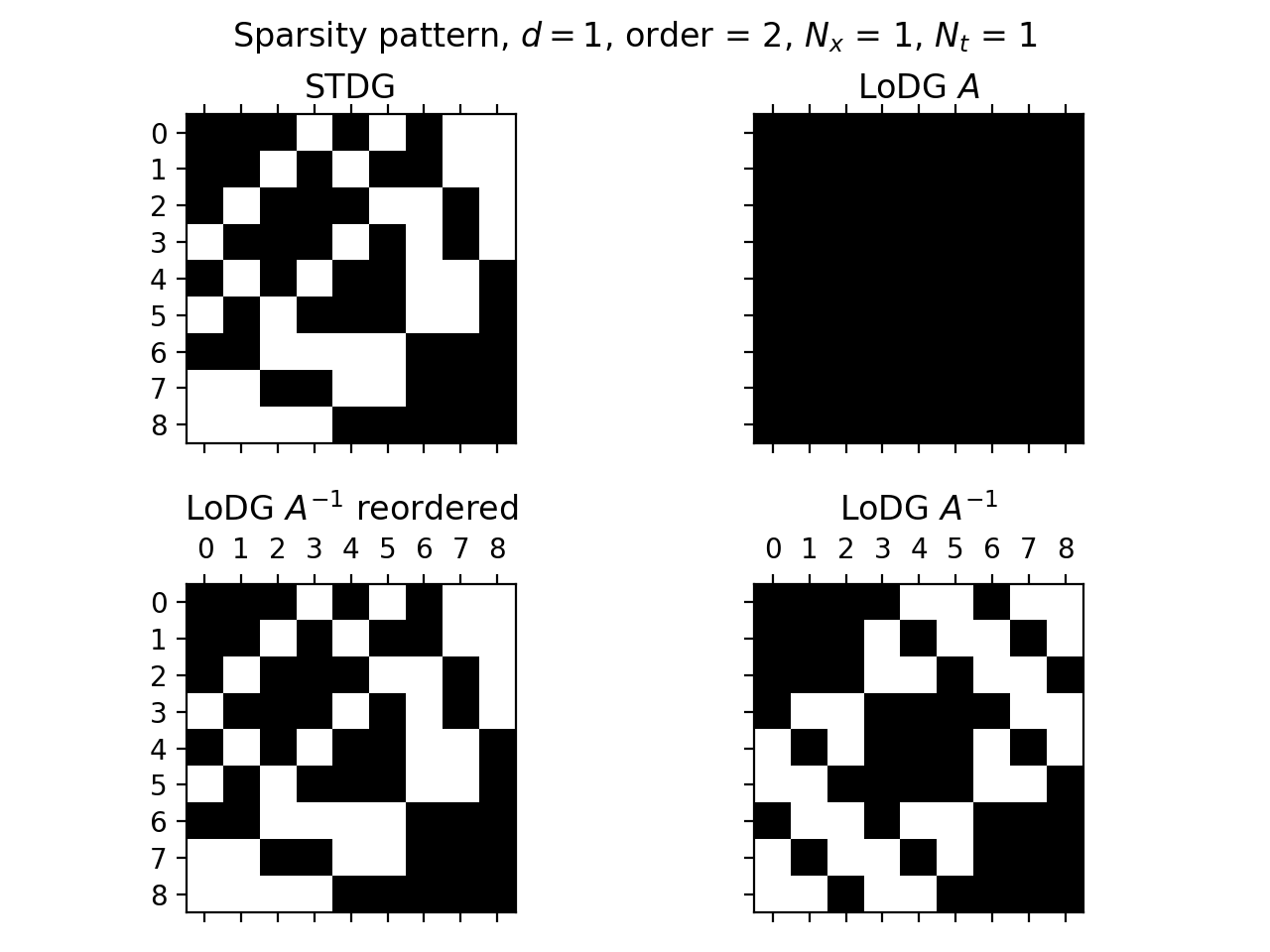}
\end{subfigure}
%
%
\begin{subfigure}[b]{.475\textwidth}
\centering
\includegraphics[width=\textwidth]{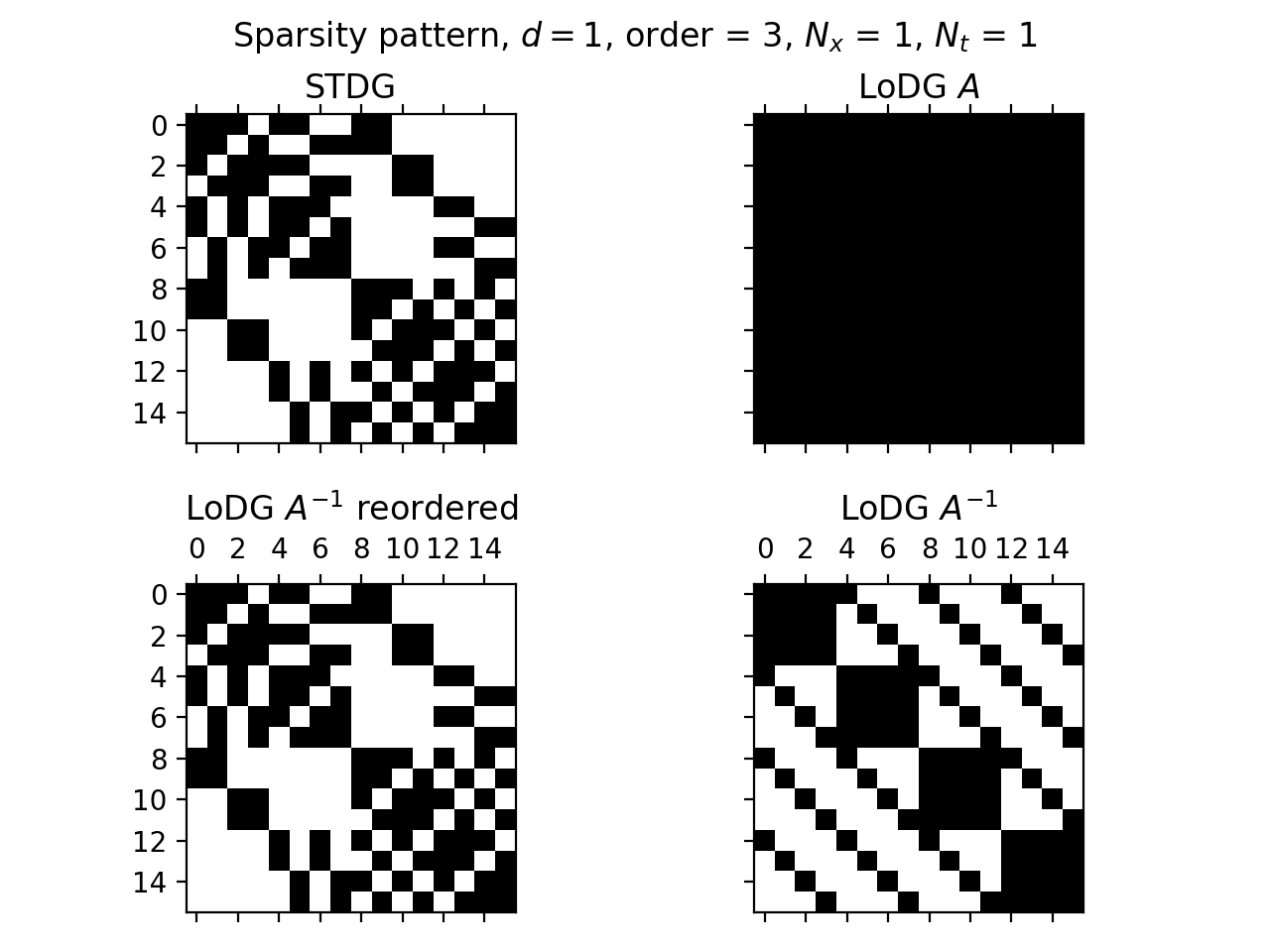}
\end{subfigure}
%
%
\begin{subfigure}[b]{.457\textwidth}
\centering
\includegraphics[width=\textwidth]{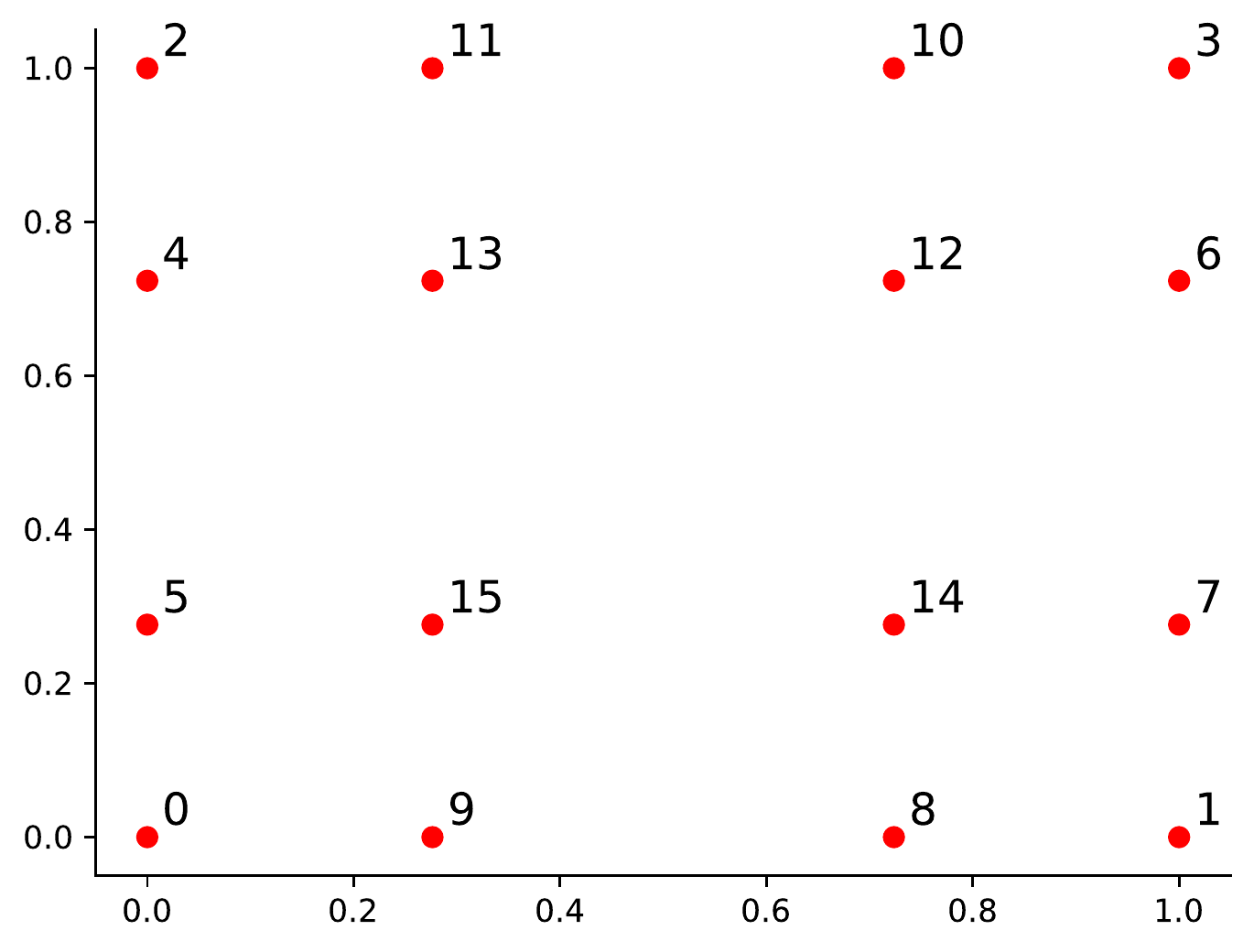}
\end{subfigure}
\caption{Sparsity patterns Jacobian of the advection problem for STDG and LoDG. 
  Node order for one space-time element ($p=3$) for STDG as generated by \dunefem (bottom right).}
\label{fig:SparsityD1}
\end{figure}

Consider the 1D linear advection equation $u_t + u_x = 0$ discretized using a single element in space and time. With STDG, the discretization is generated using \dune. With LoDG, the spatial terms are generated with \dune whereas the temporal terms are set up manually as in Section \ref{ch:theory}. Figure~\ref{fig:SparsityD1} shows the sparsity patterns of the Jacobians using order 1,2 and 3 in space and time. In each figure quadruplet, the Jacobian \eqref{eq:Jacobian_ST} of STDG is shown in the top left and the Jacobian \eqref{eq:Jacobian_A} of LoDG in the top right. In the bottom right, the alternative formulation \eqref{eq:Jacobian_Ainv} is shown.

The first thing to note is that the straightforward LoDG formulation leads to a dense discretization whereas STDG is sparse. The LoDG formulation using $\mat{A}^{-1}$ is also sparse. It has the same number of nonzero elements as STDG, although their distribution is different. The explanation for this lies in the ordering of the unknowns. With LoDG, the node order is lexicographic in the temporal direction. However, the space-time element generated by \dune is as shown in the bottom right of Figure~\ref{fig:SparsityD1}, here with $\Nt = 4$.
This ordering is the result of a generic construction of the reference
elements, which is based on a recursion over the spatial
dimension $d$ starting at the $0$-dimensional reference element, i.e. a point. 
This recursion also generates a natural ordering for the basis functions,
starting with the basis functions located at points in an element and recursively down to the basis functions 
located inside the element. A detailed description of this construction is found in \cite{dedner:12}.
\noindent With a suitable permutation of the unknowns, the sparsity pattern of LoDG using $\mat{A}^{-1}$ coincides with STDG as seen in the bottom left of each figure quadruplet.

These observations suggest that the LoDG system \eqref{eq:DG_vs_RK:RK} is more expensive to work with than the STDG system \eqref{eq:DG_vs_RK:DG}, and that the $\mat{A}^{-1}$ formulation \eqref{eq:Jacobian_Ainv} is a better choice. However, the interaction of particular preconditioners and solvers with these systems may also depend on the node ordering in ways that must be deduced through careful testing.

\section{Experiments} \label{ch:experiments}
We now perform a series of numerical tests.
We start with solving the linear test equation to validate the convergence rates of the two solvers. A two-dimensional advection-diffusion test case follows, with the purpose of highlighting slight differences in the numerical solutions and particular challenges with respect to visualizing the solutions. Finally, the two and three dimensional Euler equations of gas dynamics are solved to demonstrate that both codes are capable of handling nonlinear space-time dynamics in multiple dimensions.

As mentioned previously, the spatial parts of both LoDG and STDG are generated using \dunefem. The temporal part of LoDG is implemented in \assimulo whereas \dunefem is used for the entire space-time discretization in STDG.

\subsection{Validation of Convergence Rates}

To verify that the temporal discretizations converge as expected we perform a simple test on the linear test equation,
\begin{equation} \label{eq:linear_test_eqn}
\begin{aligned}
u_t & = -u, \quad t \in (0,1], \\
u(0) & = 4.
\end{aligned}
\end{equation}
Python's sparse linear solver is used to solve any algebraic systems arising from the discretizations. The experimental order of convergence (EOC) of the pointwise error $|\ucoeff^{n+1} - u(1)|$ is shown for LoDG and STDG in Table \ref{TablePointwise}. Here, $N$ denotes the number of time steps/time elements and $\Nt \in \{ 2,3,4 \}$. Recall from Section \ref{ch:theory} that the order of LoDG, and correspondingly the superconvergence of STDG, is $2(\Nt-1)$. 
This is indeed what we observe in Table \ref{TablePointwise}. With $\Nt = 3$, the errors are approaching machine precision when $N=2^9$, hence a drop in the convergence rate is seen in Table \ref{TableEOCLGL3}. The STDG appears to be more sensitive in this respect than LoDG. The same thing happens when $\Nt=4$ and $N=2^5$, as seen in Table \ref{TableEOCLGL4}.

\begin{table}[htb]
\caption{EOC of pointwise error for LoDG and STDG applied to the test equation \eqref{eq:linear_test_eqn}.}
\label{TablePointwise}
\renewcommand\arraystretch{1.2}
\centering
\begin{subtable}[t]{0.3\textwidth}
\centering
\caption{$\Nt = 2$}
\label{TableEOCLGL2}
\begin{tabular}{ | c | c | c |}
\hline
 N & Lobatto & DG-SEM \\ 
 \hline 
 $2^4$ & 1.93 & 1.93 \\
 $2^5$ & 1.97 & 1.97 \\
 $2^6$ & 1.98 & 1.98 \\
 $2^7$ & 1.99 & 1.99 \\
 $2^8$ & 1.99 & 1.99 \\
 $2^9$ & 1.99 & 1.99 \\
 \hline
\end{tabular}
\end{subtable}
\hfil
\begin{subtable}[t]{0.3\textwidth}
\centering
\caption{$\Nt = 3$}
\label{TableEOCLGL3}
\begin{tabular}{ | c | c | c |}
\hline
 N & Lobatto & DG-SEM \\ 
 \hline 
 $2^4$ & 3.96 & 3.96 \\
 $2^5$ & 3.98 & 3.98 \\
 $2^6$ & 3.99 & 3.99 \\
 $2^7$ & 3.99 & 4.01 \\
 $2^8$ & 3.99 & 4.28 \\
 $2^9$ & 4.05 & 0.47 \\
 \hline
\end{tabular}
\end{subtable}
\hfil
\begin{subtable}[t]{0.3\textwidth}
\centering
\caption{$\Nt = 4$}
\label{TableEOCLGL4}
\begin{tabular}{ | c | c | c |}
\hline
 N & Lobatto & DG-SEM \\ 
 \hline 
 $2^4$ & 5.98 & 5.96 \\
 $2^5$ & 6.64 & 4.22 \\
 \hline
\end{tabular}
\end{subtable}
\end{table}


We now repeat the experiment but measure the EOC via the $L_2$ norm $\| \cdot \|_{L_2[0,1]}$. This type of error measurement is straightforward to perform with the STDG code. However, the LoDG implementation does not by default save the intermediate RK stages necessary to perform the computation of the $L_2$ error. We expect the EOC to be given by $\Nt$. Table \ref{TableL2} shows that this indeed is observed. For $\Nt = 4$ with $N=2^8$ time elements the convergence rate drops due to very small errors, as seen in Table \ref{TableEOCLGL4L2}. Again, STDG appears to be more sensitive to this phenomenon than LoDG.

\begin{table}[htb]
\caption{EOC of $L_2$ error for LoDG and STDG applied to the test equation \eqref{eq:linear_test_eqn}.}
\label{TableL2}
\renewcommand\arraystretch{1.2}
\centering
\begin{subtable}[t]{0.3\textwidth}
\centering
\caption{$\Nt = 2$}
\label{TableEOCLGL2L2}
\begin{tabular}{ | c | c | c |}
\hline
 N & Lobatto & DG-SEM \\ 
 \hline 
 $2^4$ & 1.96 & 1.96 \\
 $2^5$ & 1.98 & 1.98 \\
 $2^6$ & 1.99 & 1.99 \\
 $2^7$ & 2.0  & 2.0 \\
 $2^8$ & 2.0  & 2.0 \\
 $2^9$ & 2.0  & 2.0 \\
 \hline
\end{tabular}
\end{subtable}
\hfil
\begin{subtable}[t]{0.3\textwidth}
\centering
\caption{$\Nt = 3$}
\label{TableEOCLGL3L2}
\begin{tabular}{ | c | c | c |}
\hline
 N & Lobatto & DG-SEM \\ 
 \hline 
 $2^4$ & 2.98 & 2.98 \\
 $2^5$ & 2.99 & 2.99 \\
 $2^6$ & 2.99 & 2.99 \\
 $2^7$ & 3.0  & 3.0 \\
 $2^8$ & 3.0  & 3.0 \\
 $2^9$ & 3.0  & 3.0 \\
 \hline
\end{tabular}
\end{subtable}
\hfil
\begin{subtable}[t]{0.3\textwidth}
\centering
\caption{$\Nt = 4$}
\label{TableEOCLGL4L2}
\begin{tabular}{ | c | c | c |}
\hline
 N & Lobatto & DG-SEM \\ 
 \hline 
 $2^4$ & 3.99 & 3.99 \\
 $2^5$ & 4.0  & 4.0 \\
 $2^6$ & 4.0  & 4.0 \\
 $2^7$ & 4.0  & 4.0 \\
 $2^8$ & 4.0  & 3.6 \\
 \hline
\end{tabular}
\end{subtable}
\end{table}

%

\subsection{Advection-Diffusion}

The next test case is the linear advection-diffusion problem in two dimensions;
\begin{equation} \label{advdiffeq}
\begin{alignedat}{3}
\partial_t \vecu + \mathbf{b} \cdot \nabla \vecu - \varepsilon \Delta \vecu &= 0 \qquad && \mbox{in } ( \Omega \times (0,T]),~ \Omega \subset \R^2, \\
\vecu(0) &= \vecu_0 \qquad && \mbox{in } \Omega.
\end{alignedat}
\end{equation}
We test both implementations for the rotating pulse problem with analytic solution
\begin{align*}
\vecu(t,\mathbf{x}) & = \frac{0.004}{0.004+4\varepsilon t} \exp\left( -\frac{x_q^2 + y_q^2}{0.004 + 4 \varepsilon t} \right),\\
x_q & = x_0\cos(4t) + y_0\sin(4t) + 0.25,\\
y_q & = - x_0\sin(4t) + y_0\cos(4t).
\end{align*}
Here, $x_0 = x - 0.5$, $y_0 = y - 0.5$, $\mathbf{b} = [-4y_0,4x_0]$, $\varepsilon = 0.001$, and $(t,\mathbf{x}) \in [0,1] \times [0,1]^2$. The initial condition is given by $\vecu(0,\mathbf{x})$ and we apply periodic boundary conditions in space. The linear systems arising from the discretizations are solved using built-in routines in \dune and \assimulo. Thus, despite the mathematical equivalence of LoDG and STDG, we do not expect the two codes to yield identical solutions.

\begin{figure}[h!]
\centering
\begin{subfigure}[b]{.32\textwidth}
\centering
\includegraphics[width=\textwidth]{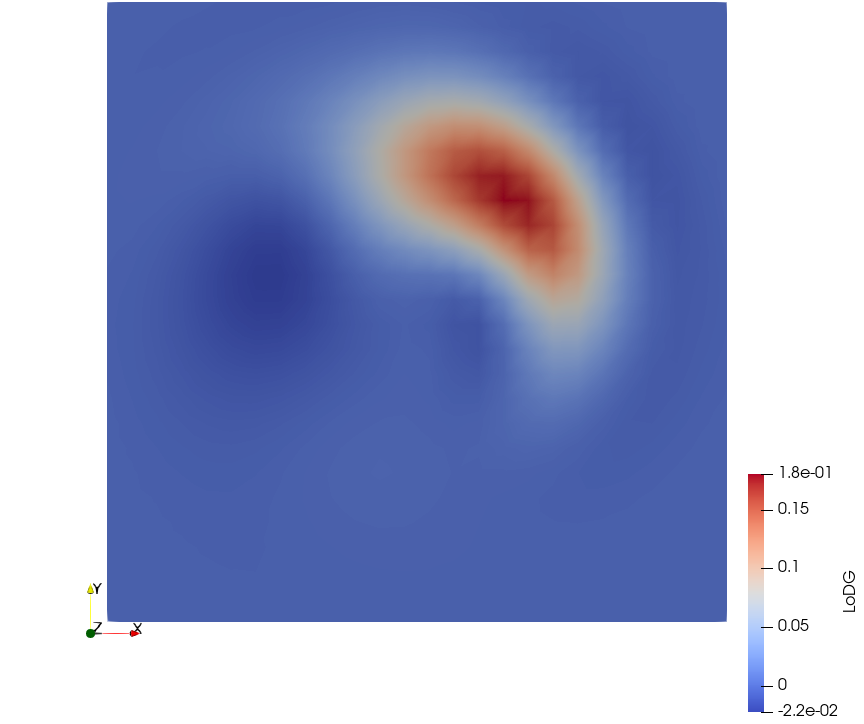}
\caption{LoDG, $\Nt=2$}
\end{subfigure}
%
%
\begin{subfigure}[b]{.32\textwidth}
\centering
\includegraphics[width=\textwidth]{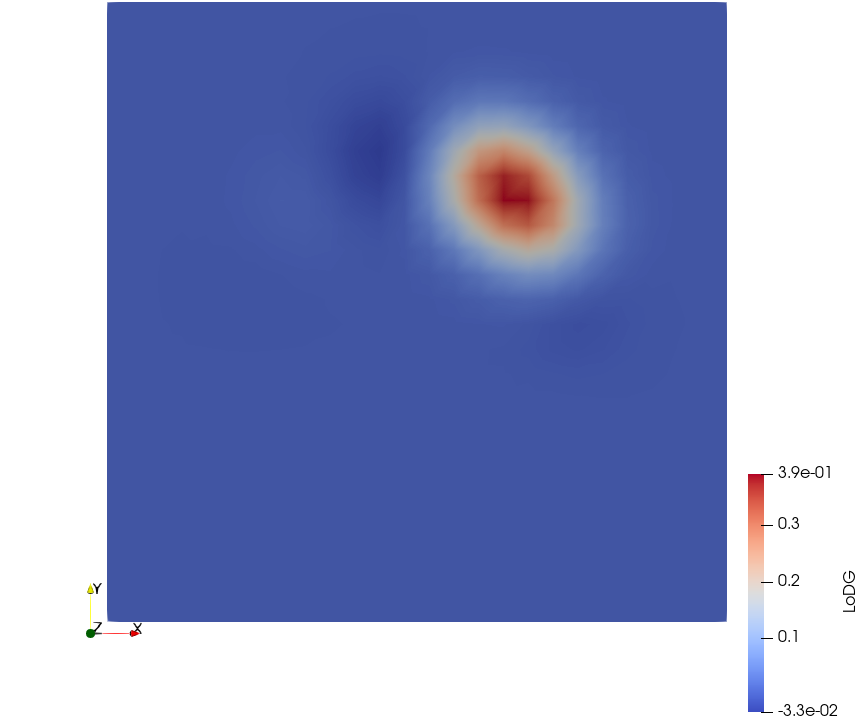}
\caption{LoDG, $\Nt=3$}
\end{subfigure}
%
%
\begin{subfigure}[b]{.32\textwidth}
\centering
\includegraphics[width=\textwidth]{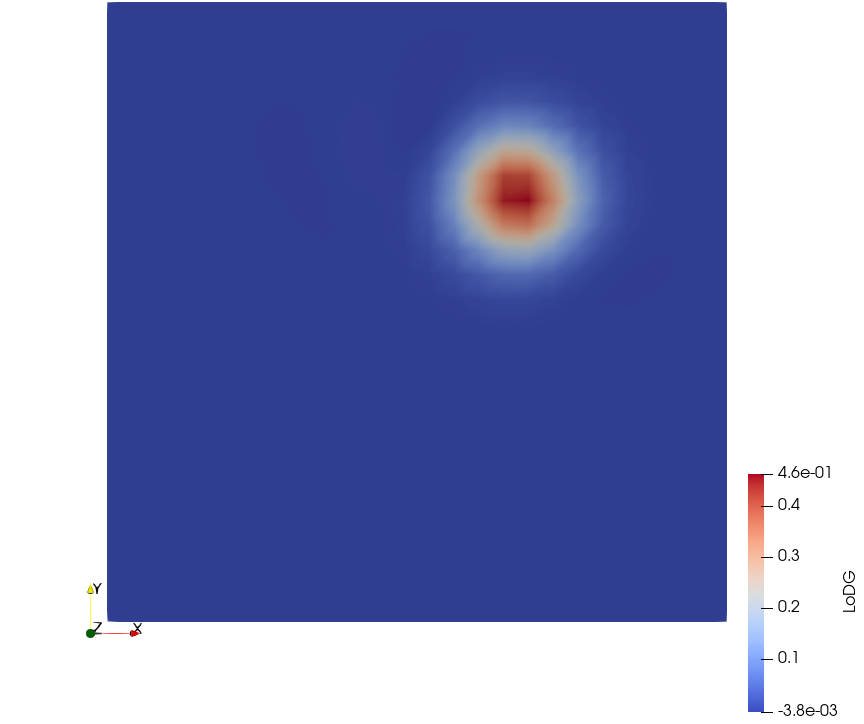}
\caption{LoDG, $\Nt=4$}
\end{subfigure}
\begin{subfigure}[b]{.32\textwidth}
\centering
\includegraphics[width=\textwidth]{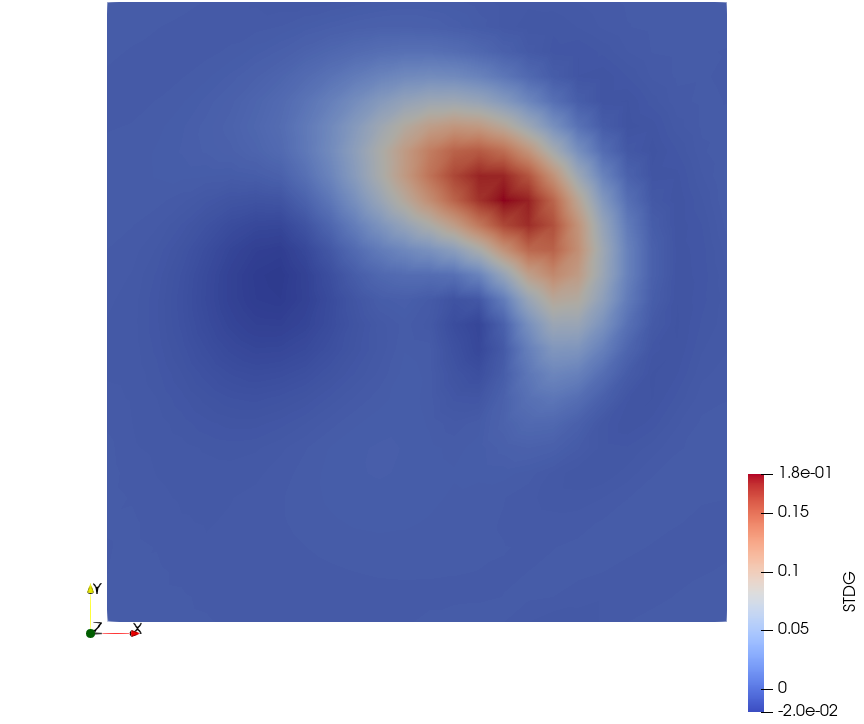}
\caption{STDG, $\Nt=2$}
\end{subfigure}
%
%
\begin{subfigure}[b]{.32\textwidth}
\centering
\includegraphics[width=\textwidth]{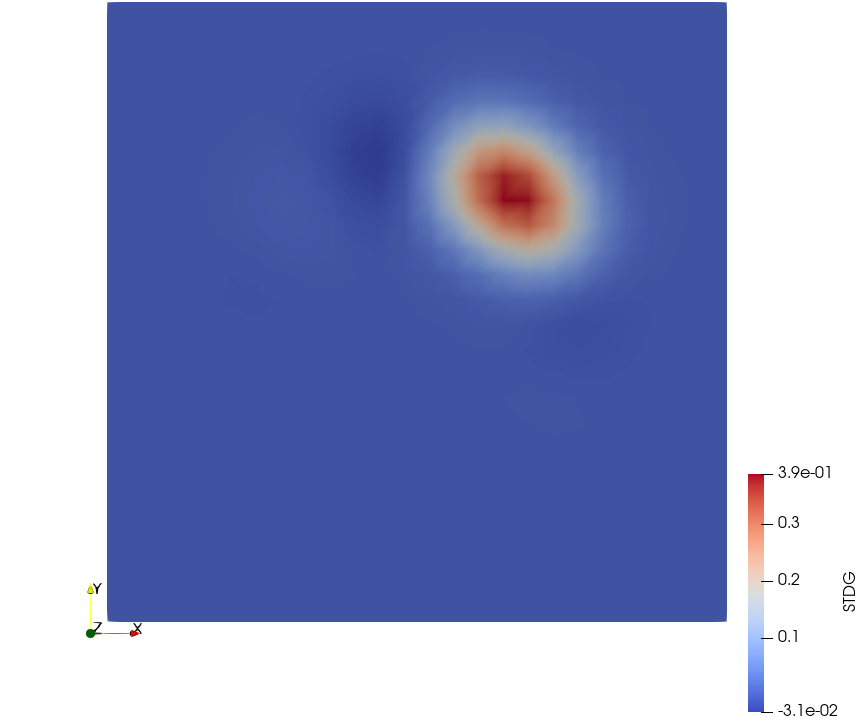}
\caption{STDG, $\Nt=3$}
\end{subfigure}
%
%
\begin{subfigure}[b]{.32\textwidth}
\centering
\includegraphics[width=\textwidth]{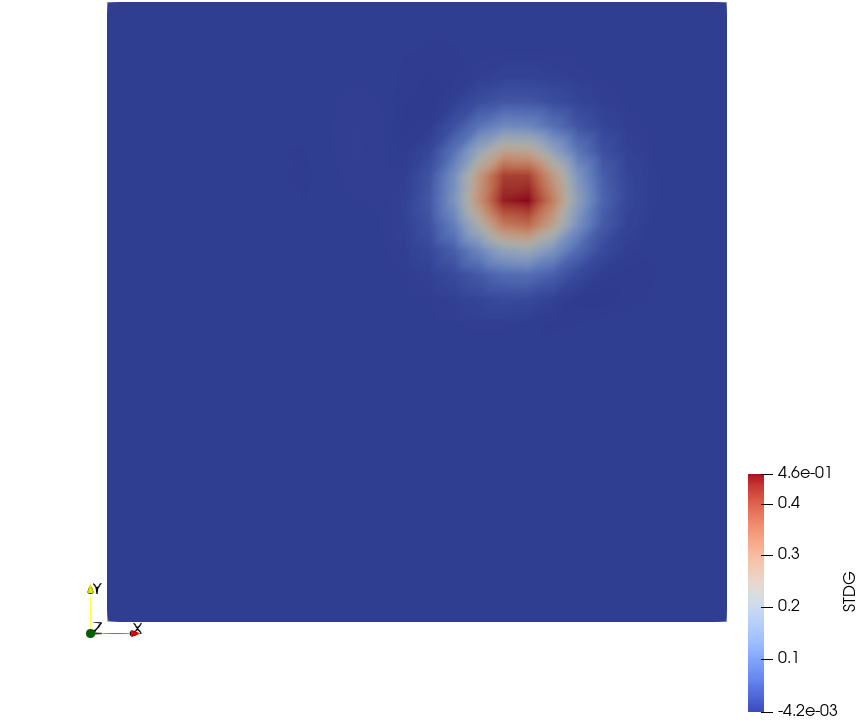}
\caption{STDG, $\Nt=4$}
\end{subfigure}
\caption{Numerical solution of a rotating pulse subject to the advection-diffusion equation \eqref{advdiffeq} using LoDG (top) and STDG (bottom).}
\label{fig:AdvDiff}
\end{figure}

The numerical solutions obtained by the two codes with $\Nt \in \{2,3,4\}$ are shown in Figure~\ref{fig:AdvDiff}. Here, a uniform mesh is used in space with $\Delta x = \Delta y = 0.04$. Time steps of uniform size $\Delta t = 0.1$ are used throughout the simulation. With $\Nt=2$ the problem is significantly under-resolved, leading to a smeared solution. As $\Nt$ is increased, this phenomenon is reduced. To the eye, the numerical solutions using the two codes are barely distinguishable.

To get a more detailed comparison of the numerical results obtained by the two implementations we compare their spatial $L_2$ error in the final time point, $t=1$. This time we vary the space-time grid with $\Delta x = \Delta y = \Delta t = 1/N$. The errors and the EOC are shown in Table~\ref{TableL2Errors}. Notice that the $L_2$ errors are very similar, although not identical, testifying to the influence of the different solvers of the algebraic equations. Note also that the behaviour of the EOC is less clear than it was for the linear test equation. In this experiment we have refined space and time simultaneously, and therefore do not have a theoretical convergence result to rely on. The results indicate a convergence rate higher than $\Nt$, although not quite as high as $2(\Nt-1)$.

\begin{table}[ht!]
\renewcommand\arraystretch{1.2}
\centering
\caption{Errors and EOC for the advection-diffusion problem \eqref{advdiffeq}.}
\label{TableL2Errors}
\begin{subtable}{0.8\textwidth}
\centering
\begin{tabular}{ | >{\centering}p{1cm} | >{\centering}p{1.5cm} | >{\centering}p{1.5cm} | >{\centering}p{1.5cm} | >{\centering}p{1.5cm} | >{\centering}p{1.5cm} | >{\centering\arraybackslash}p{1.5cm} |}
\hline
& \multicolumn{2}{c|}{$\Nt=2$} & \multicolumn{2}{c|}{$\Nt=3$} & \multicolumn{2}{c|}{$\Nt=4$} \\
\hline
 $N$ & LoDG & STDG & LoDG & STDG & LoDG & STDG \\ 
 \hline 
 $2^2$ & 8.94E-2 & 7.28E-2 & 4.45E-2 & 4.37E-2 & 2.68E-2 & 2.69E-2 \\ 
 \hline 
 $2^3$ & 4.66E-2 & 4.46E-2 & 2.42E-2 & 2.41E-2 & 6.05E-3 & 6.04E-3 \\ 
 \hline 
 $2^4$ & 3.49E-2 & 3.39E-2 & 5.36E-3 & 5.38E-3 & 4.92E-4 & 4.93E-4 \\ 
 \hline 
 $2^5$ & 1.86E-2 & 1.84E-2 & 5.85E-4 & 5.94E-4 & 1.06E-5 & 9.88E-6 \\   
 \hline 
\end{tabular}
\caption{Error}
\end{subtable}

\medskip
\begin{subtable}{0.8\textwidth}
\centering
\begin{tabular}{ | >{\centering}p{1cm} | >{\centering}p{1.5cm} | >{\centering}p{1.5cm} | >{\centering}p{1.5cm} | >{\centering}p{1.5cm} | >{\centering}p{1.5cm} | >{\centering\arraybackslash}p{1.5cm} |}
\hline
& \multicolumn{2}{c|}{$\Nt=2$} & \multicolumn{2}{c|}{$\Nt=3$} & \multicolumn{2}{c|}{$\Nt=4$} \\
\hline
 $N$ & LoDG & STDG & LoDG & STDG & LoDG & STDG \\ 
 \hline 
 $2^3$ & 0.9 & 0.7 & 0.9 & 0.9 & 2.1 & 2.2 \\ 
 \hline 
 $2^4$ & 0.4 & 0.4 & 2.2 & 2.2 & 3.6 & 3.6 \\ 
 \hline 
 $2^5$ & 0.9 & 0.9 & 3.2 & 3.2 & 5.5 & 5.6 \\   
 \hline 
\end{tabular}
\caption{EOC}
\end{subtable}
\end{table}

\begin{figure}[h!]
\centering
\begin{subfigure}[b]{.49\textwidth}
\centering
\includegraphics[width=\textwidth]{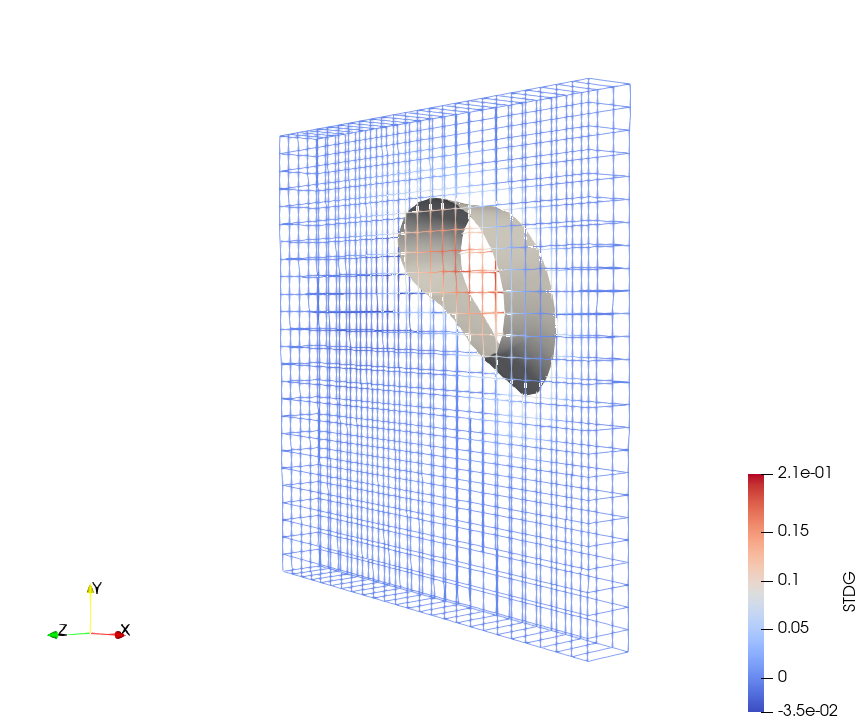}
\caption{$\Nt=2$}\label{fig:AdvDiffSiceO1}
\end{subfigure}
%
%
\begin{subfigure}[b]{.49\textwidth}
\centering
\includegraphics[width=\textwidth]{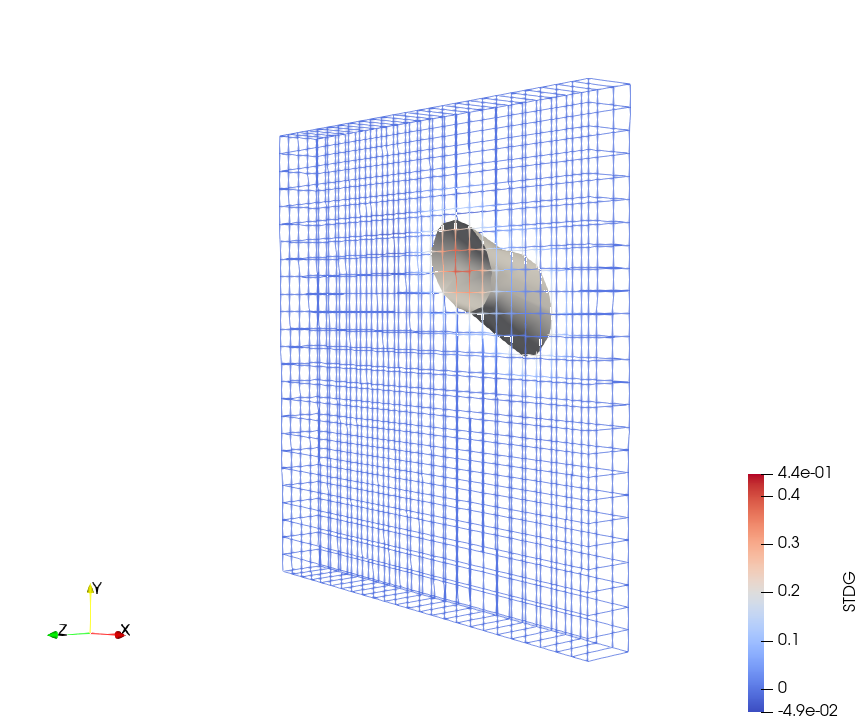}
\caption{$\Nt=3$}\label{fig:AdvDiffSiceO2}
\end{subfigure}
%
%
\begin{subfigure}[b]{.49\textwidth}
\centering
\includegraphics[width=\textwidth]{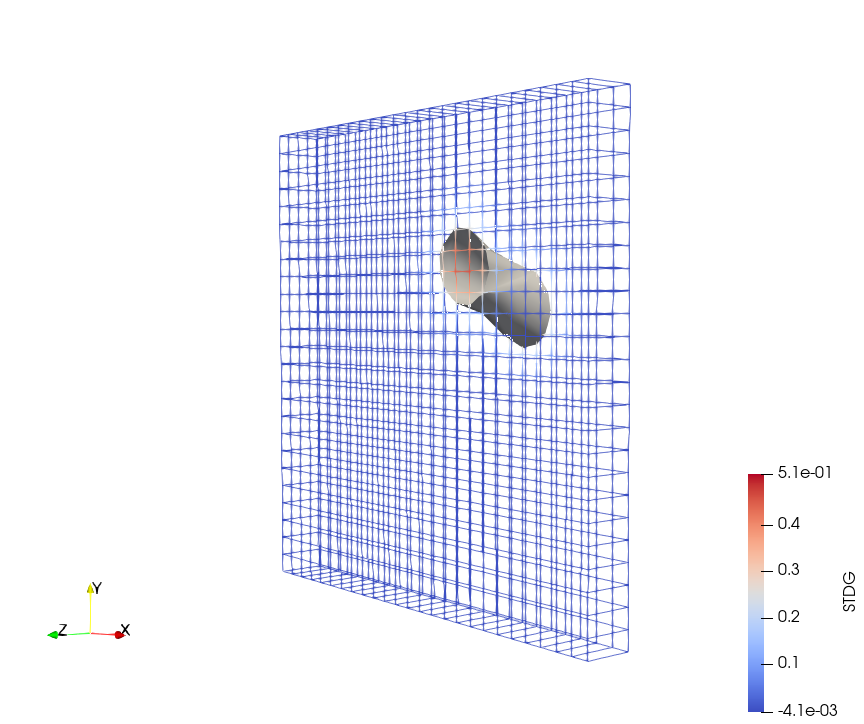} %
\caption{$\Nt=4$}\label{fig:AdvDiffSiceO3}
\end{subfigure}
%
%
\begin{subfigure}[b]{.49\textwidth}
\centering
\includegraphics[width=\textwidth]{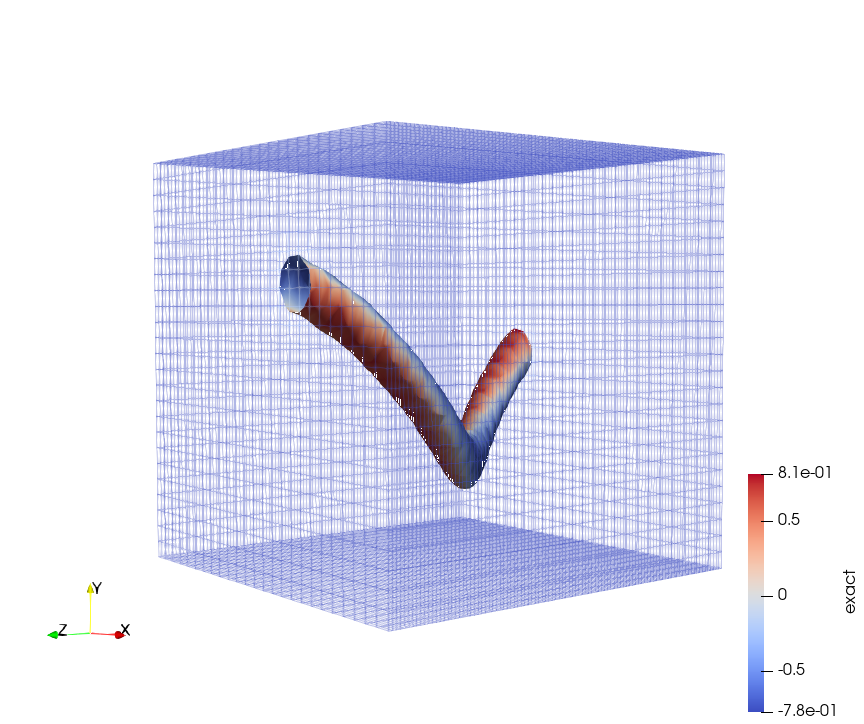}
\caption{Exact solution in space-time}\label{ig:AdvDiff3dExact}
\end{subfigure}
\caption{3D space-time visualization of the rotating pulse produced with the STDG code.}
\label{fig:AdvDiffSlice}
\end{figure}

Finally, we highlight a feature of the STDG code that may be of use in certain situations. Since this code returns all points in a given time element (or equivalently, all intermediate RK stages in each time step), these can be visualized using a 3D plotting software, thereby obtaining a space-time visualization of the solution. These stages are usually discarded by ODE solvers for efficiency reasons. The visualization is done for a single time step in Figure~\ref{fig:AdvDiffSlice}. Here, the exact solution over the whole space-time domain is also shown for reference.


\subsection{Euler Equations}

A prime example for evolution equations are the \textit{Euler equations
of gas dynamics}. They are derived from the conservation of mass,
momentum, and energy of a compressible inviscid fluid.
In Eulerian coordinates they have the form:
\begin{equation}  \label{eulereq}
\begin{alignedat}{3}
\partial_t \vecu + \nabla \cdot \Fc(\vecu) &= 0 \qquad && \mbox{in } (\Omega \times (0,T]),~\Omega \subset \R^d, \quad d \in \{1,2,3\}, \\
\vecu(0) &= \vecu_0 \qquad && \mbox{in } \Omega, 
\end{alignedat}
\end{equation}
where the vector of the conservative variables has the form
\begin{eqnarray}\label{defu}
\vecu = \left ( \begin{array}{c} 
            \rho \\
            \rho \vecv \\
            \energy 
           \end{array} \right ), 
       \quad \rho \vecv = (\rho v_1,\ldots,\rho v_d)^T, \quad \energy = \rho \mathcal{E},
\end{eqnarray}
augmented with suitable boundary conditions (which are discussed in detail in \cite{Birken2021}).
Here, $\rho$ denotes the density of the fluid, $\vecv$ the velocity, $\energy$ the internal energy, and $\mathcal{E}$ the total energy.
%
%
The convective flux function $\Fc(\vecu) := (\vect{f}_1(\vecu),\ldots,\vect{f}_d(\vecu))$
has for $i=1,...,d$ the form
\begin{eqnarray*}
\vect{f}_i(\vecu) := \left ( \begin{array}{c} 
             \vecu_{i+1} \\
             \vecu_{i+1} \vecu_2 / \vecu_1 \, + \delta_{i,1}\, \pressure(\vecu) \\
             \vdots \\
             \vecu_{i+1} \vecu_{d+1} / \vecu_1 \, + \delta_{i,d}\, \pressure(\vecu) \\
             (\vecu_{d+2} + \pressure(\vecu)) \, \vecu_{i+1}/\vecu_1 
           \end{array} \right ), 
\end{eqnarray*}
where $\delta_{i,j}$ is the Kronecker delta.
For example, choosing $d=3$ and directly using $\vecu$ from (\ref{defu})
we obtain the three flux functions
\begin{eqnarray*}
\vect{f}_1(\vecu) = \left ( \begin{array}{c} 
             \rho v_1 \\
             \rho v_1^2 + \pressure \\
             \rho v_1 v_2  \\
             \rho v_1 v_3  \\
             (\energy + \pressure) v_1
           \end{array} \right ), 
           \quad 
\vect{f}_2(\vecu) = \left ( \begin{array}{c} 
             \rho v_2 \\
             \rho v_2 v_1  \\
             \rho v_2^2 + \pressure \\
             \rho v_2 v_3  \\
             (\energy + \pressure) v_2
           \end{array} \right ),
           \quad 
\vect{f}_3(\vecu) = \left ( \begin{array}{c} 
             \rho v_3 \\
             \rho v_3 v_1  \\
             \rho v_3 v_2  \\
             \rho v_3^2 + \pressure \\
             (\energy + \pressure) v_3
           \end{array} \right ).
\end{eqnarray*}

We first consider the two-dimensional Euler equations 
with periodic boundary conditions and vortex initial condition
\begin{align*}
\rho & = \left(1 - S^2(\gamma-1)M^2\frac{\exp(f)}{(8\pi^2)}\right)^{\frac{1}{\gamma-1}} ,\\
v_1 & = 1 - S \mathbf{x}_1 \frac{\exp\left(\frac{f}{2}\right)}{2\pi}, \quad v_2 = S \mathbf{x}_0 \frac{\exp\left(\frac{f}{2}\right)}{2\pi}, \\
\energy & = \frac{\pressure}{\gamma-1} + 0.5\frac{v_1^2+v_2^2}{\rho},~ \pressure = \frac{\rho^{\gamma}}{\gamma M^2},
\end{align*}
with vortex strength $S = 5$, Mach number $M=0.5$, $\gamma=1.4$ and $f = 1 - \mathbf{x}_0^2 - \mathbf{x}_1^2$.

\begin{figure}[h!]
\centering
\begin{subfigure}[b]{.32\textwidth}
\centering
\includegraphics[width=\textwidth]{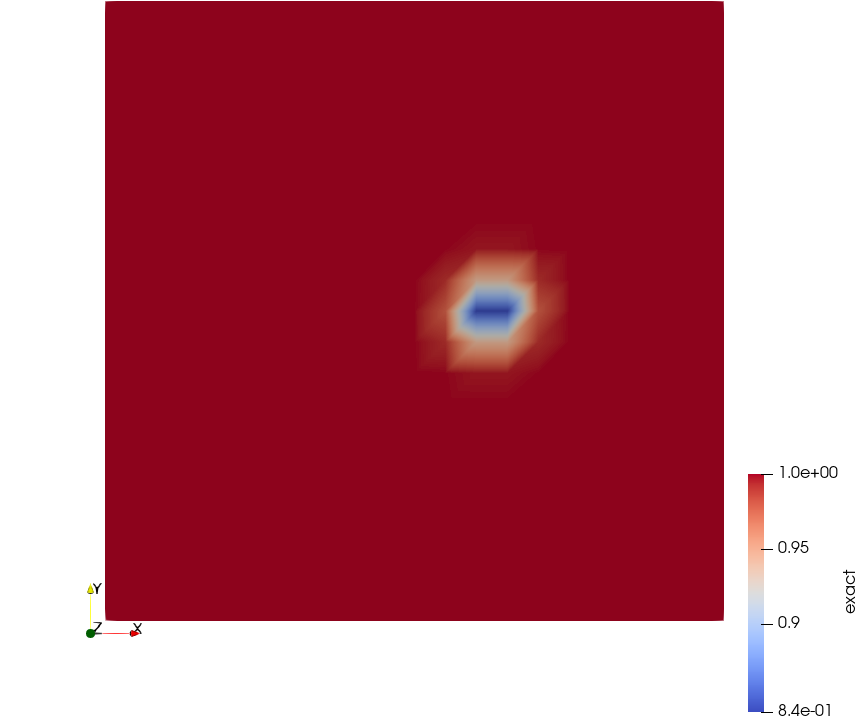}
\caption{Exact, $\Nt=2$}\label{fig:EulerO1Exact}
\end{subfigure}
%
%
\begin{subfigure}[b]{.32\textwidth}
\centering
\includegraphics[width=\textwidth]{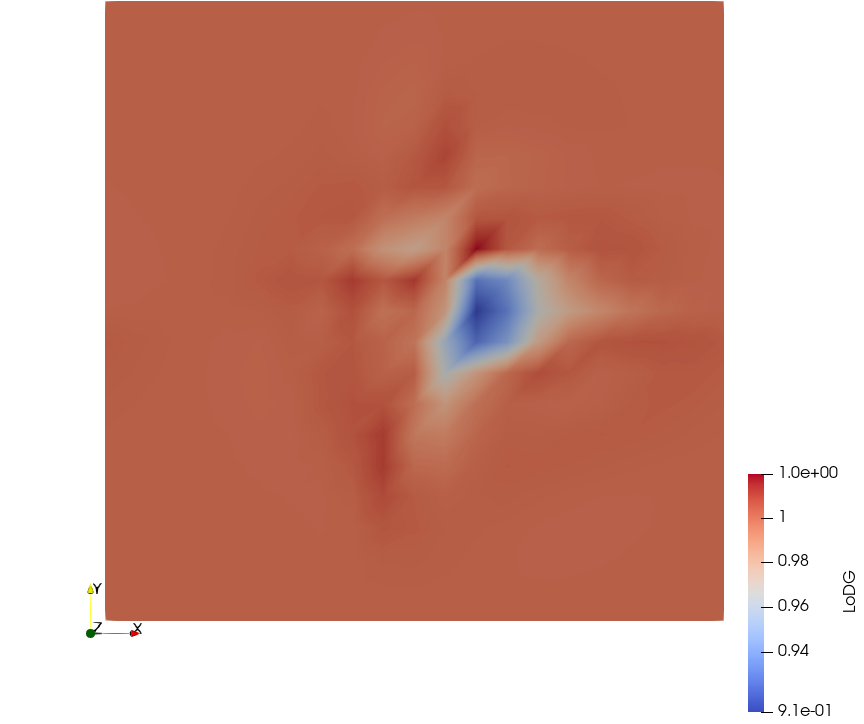}
\caption{LoDG, $\Nt=2$}\label{fig:EulerO1LoDG}
\end{subfigure}
%
%
\begin{subfigure}[b]{.32\textwidth}
\centering
\includegraphics[width=\textwidth]{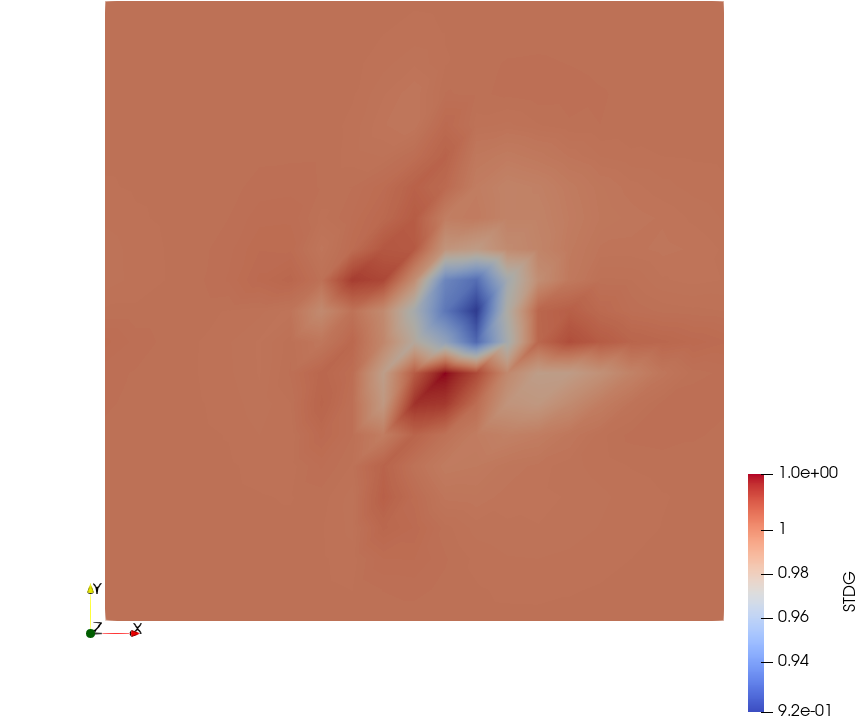}
\caption{STDG, $\Nt=2$}\label{fig:EulerO1STDG}
\end{subfigure}
\begin{subfigure}[b]{.32\textwidth}
\centering
\includegraphics[width=\textwidth]{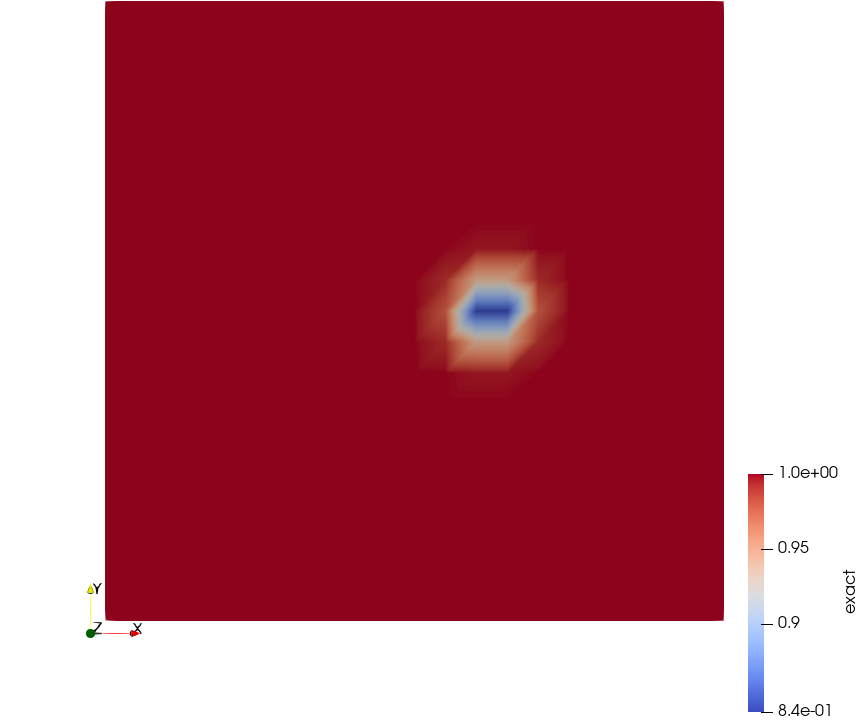}
\caption{Exact, $\Nt=3$}\label{fig:EulerO2Exact}
\end{subfigure}
%
%
\begin{subfigure}[b]{.32\textwidth}
\centering
\includegraphics[width=\textwidth]{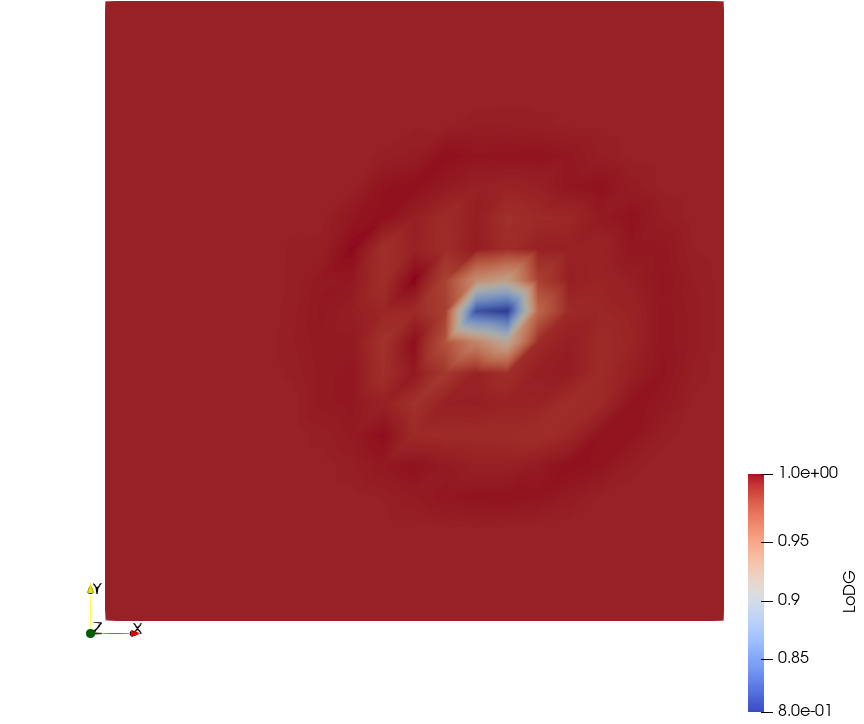}
\caption{LoDG, $\Nt=3$}\label{fig:EulerO2LoDG}
\end{subfigure}
%
%
\begin{subfigure}[b]{.32\textwidth}
\centering
\includegraphics[width=\textwidth]{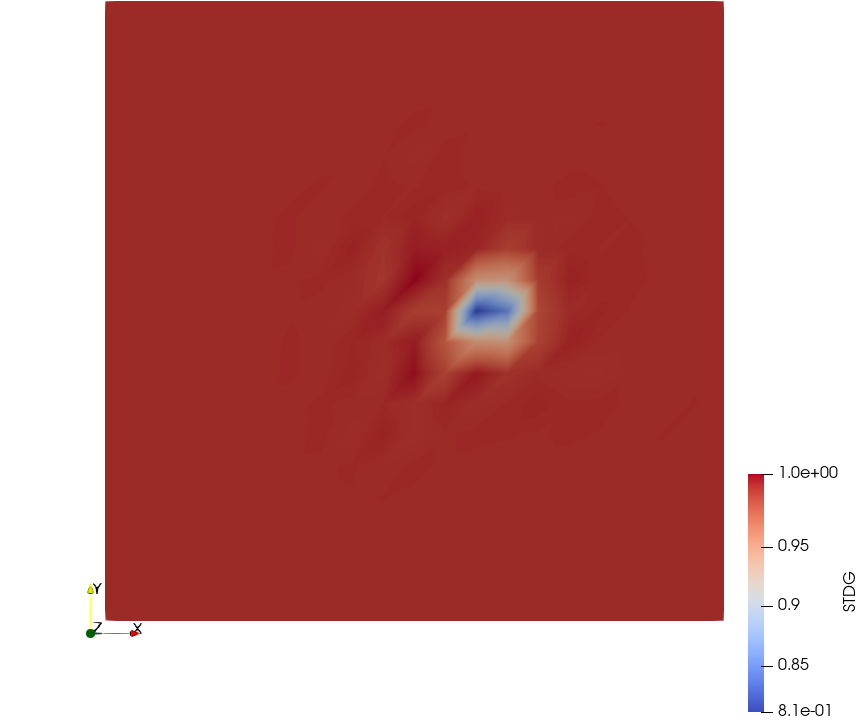}
\caption{STDG, $\Nt=3$}\label{fig:EulerO2STDG}
\end{subfigure}
\caption{Numerical solution of $\rho$ for a vortex problem subject to the 2D Euler equations \eqref{eulereq} using LoDG and STDG.}
\label{fig:Euler2D}
\end{figure}

The numerical solutions obtained by the two codes with $\Nt \in \{2,3\}$ are shown in Figure~\ref{fig:Euler2D}. Here, a uniform mesh on the space-time domain $[-10,10]^2 \times (0,2.5]$ is used in space with $\Delta x = \Delta y = 0.04$. Time steps of uniform size $\Delta t = 0.01$ are used throughout the simulation. Again, the problem is significantly under-resolved with $\Nt=2$, leading to a smeared solution. As $\Nt$ is increased, this phenomenon is reduced. Some differences can be seen in the numerical results for the two implementations. This is likely caused by the fact that two different solvers, inherent to \assimulo and \dune respectively, are used for the nonlinear systems arising from the discretizations. Due to the differences between these solvers we can in general not expect identical numerical solutions despite the mathematical equivalence of the two algorithms.

%

To show the potential of our STDG code we present a 3D Euler test case with $\Nt = 3$. We use periodic boundary conditions and a smooth bubble advection
initial condition
\begin{align*}
\rho & = 
\begin{cases}
0.5,\quad & \hat{x} > 1.0,\\
0.25 (\cos(\hat{x} \pi) +1)^2 + 0.5, \quad & \hat{x} \le 1.0,
\end{cases}\\
v_1 & = \cos\bigl(\frac{\pi}{5}\bigr),\quad v_2 = \sin\bigl(\frac{\pi}{5}\bigr), \quad v_3 = \sin\bigl(\frac{\pi}{5}\bigr) \\
\energy & = \frac{p}{\gamma-1} + 0.5\rho (v_1^2 + v_2^2 + v_3^2),
\end{align*}
with $p = 0.3$, $\gamma=1.4$ and $\hat{x} = 16\sum_{i=1}^{d} (\mathbf{x}_i - 0.25 - t v_i)^2$.

\begin{figure}[h!]
\centering
\begin{subfigure}[b]{.49\textwidth}
\centering
\includegraphics[width=\textwidth]{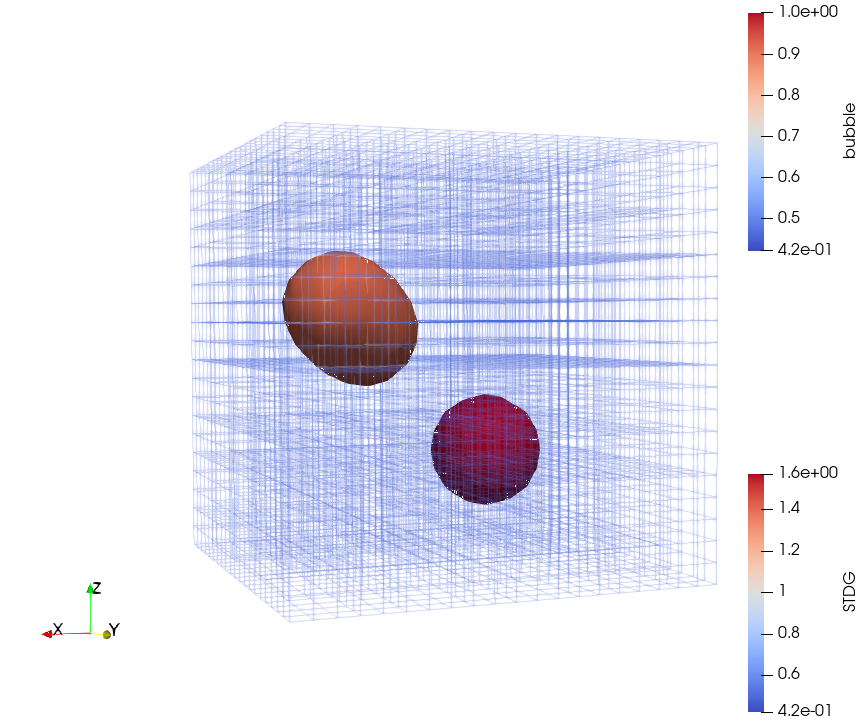}
\caption{$\Nt=3$, $\Delta t = 0.2$}\label{fig:EulerO2ExactrhoNew}
\end{subfigure}
%
%
\begin{subfigure}[b]{.49\textwidth}
\centering
\includegraphics[width=\textwidth]{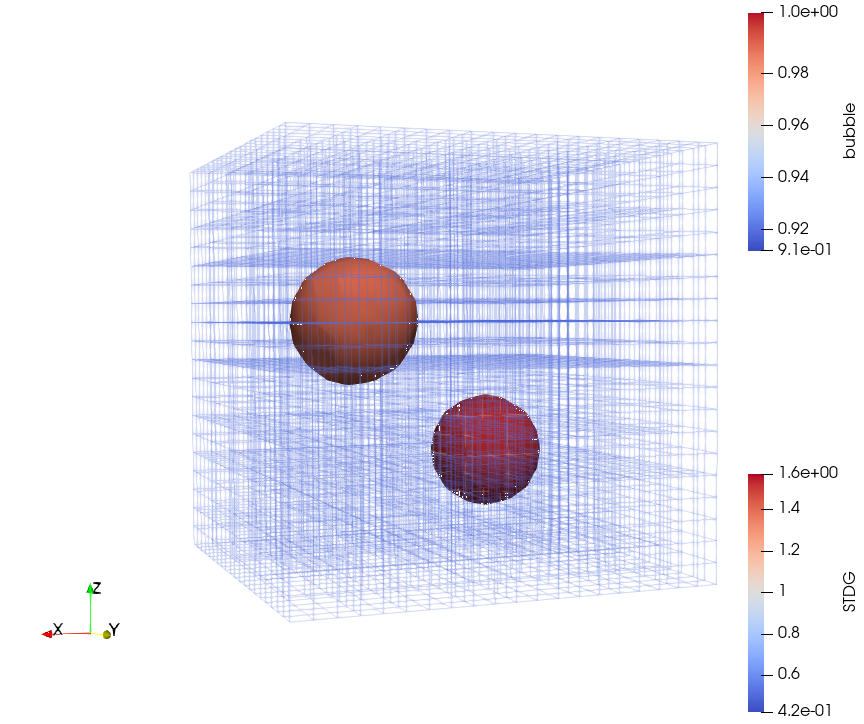}
\caption{$\Nt=3$, $\Delta t = 0.05$}\label{fig:EulerO2STDGrhoNew}
\end{subfigure}
\caption{Numerical solution of $\rho$ for a smooth bubble advection problem subject to the 3D Euler equations \eqref{eulereq} using STDG, initial condition (lower right bubble) and numerical solution (upper left bubble).}
\label{fig:Euler3DO2rhoNew}
\end{figure}

We consider the space-time domain $[0,1]^3 \times (0,0.6]$ and use a uniform mesh in space with $\Delta x = \Delta y = \Delta z = 0.05$ and time steps of uniform size $\Delta t = 0.2$ and $\Delta t = 0.05$ throughout the simulation.
The results have been computed on the LUNARC Aurora cluster at Lund University using 
$640$ Intel Xeon E5-2650 v3 processors and a Newton-GRMES solver with SOR preconditioning 
based on the PETSc library \cite{petsc-user-ref, petsc-web-page}. 
On average we observe $3$ Newton iterations and about $25$ linear iterations per timestep for $\Delta t = 0.05$ and about $50$ linear iterations per timestep for $\Delta t = 0.2$. 
The initial condition and the final time element of the density can be seen in 
Figure~\ref{fig:Euler3DO2rhoNew}. They show that some care must be taken for the numerical solution to be properly resolved in time, with a time step size $\Delta t = 0.2$ resulting in a smeared out solution with an $L_2$ error of $0.05$ in the last time point. When decreasing the time step to $\Delta t = 0.05$ the numerical solution is more accurate with an $L_2$ error of $0.007$. These results
show the potential of the \dune code even for 4D problems. Recently, 4D problems have been taken into consideration \cite{Frontin2021}, but to the best of our knowledge this is the first 4D DG-SEM implementation available publicly.

\section{Conclusions} \label{ch:conclusions}

In this paper we have presented a comparison of the theoretical and practical aspects of two different space-time DG-SEM implementations. DG-SEM in time using an upwind numerical flux is equivalent to the Lobatto IIIC family of Runge-Kutta methods in the sense that they yield the same numerical solution when solved exactly. The two methods consequently have identical order, stability and convergence properties. However, since they emanate from different research communities, the terminology used to describe them are different. In this article we have made an effort to bridge this gap and to clarify the mathematical connections between the two methods.

Having two equivalent formulations of the same discretization raises the question of which path to take towards its implementation. Two strategies are immediately obvious: Either use the method of lines with DG-SEM in space and Lobatto IIIC in time (LoDG) or use a space-time DG-SEM approach (STDG). LoDG and STDG have been implemented using \dunefem and \assimulo and the codes are available as supplementary material to this paper; see \nameref{Appendix_Installation} for details.

Despite the mathematical equivalence, there are important differences between the approaches. Not only are they described by different terminology in the literature. Additionally, the approaches lead to different systems of linear or nonlinear algebraic equations. When approximating their solutions using iterative methods, the solvers interact with these systems in different ways that are difficult to predict.

On the practical side, the two approaches lead to very different software structures that lend themselves to an assortment of opportunities and challenges. An overview of algorithmic capabilities has been given, that will be useful in different simulation contexts. These include adaptive time-stepping, adaptive mesh refinement, shock capturing, preconditioning, and other computational techniques. Some of these are likely to be more readily available in one implementation than the other, especially when reusing pre-existing code. The choice of an appropriate implementation thus depends on the needs of the user as well as on the software already available.


Several general conclusions about the two implementations can be drawn: For an STDG-type implementation, the code must be able to handle stationary 4D problems. It is desirable that the temporal dimension can be discretized with a different order of accuracy than the spatial dimensions. Further, adaptive time-stepping is nontrivial. On the other hand, this approach opens the door to parallelization in the time domain in ways that are usually not available. For problems in fewer than three spatial dimensions, software packages are already available and likely highly optimized.

For an LoDG-type implementation, black-box time-stepping routines are readily availble that implement Lobatto IIIC. Adaptive time-stepping is thus no issue, and the solution can easily be obtained at any desired intermediate times. There is no extension to 4D necessary and the order of accuracy in time can be set independently of the spatial discretization. On the other hand, this approach requires the coupling of two different codes. Parallelization in time also becomes less flexible.


\appendix

\section*{Appendix A}\label{Appendix_Installation}
There exist different ways to install \dune and \assimulo. Here, we only describe the simplest 
and most straight forward way to install both, \dune and \assimulo, which is to
use a conda environment. Then the installation is done in the following way:\\
\begin{python}
conda create -n duneproject # create a new conda environment

conda activate duneproject # activate the conda environment

conda install -c conda-forge assimulo # install assimulo

pip install -U dune-fem-dg # install dune using pip, no conda package yet

conda install -c conda-forge scipy # install scipy

git clone https://bitbucket.org/nate-sime/dolfin_dg.git _dg.git # install dolfin\_dg using pip, no conda package yet
cd dolfin_dg
python3 setup.py install
\end{python}
The space-time DG-SEM code is available online at \url{https://gitlab.maths.lth.se/dune/spacetimelobattocode}.


\section*{Appendix B}\label{Appendix_Butcher}
The Butcher tableaus for the $N_{\tau}$-stage Lobatto IIIC methods with $N_{\tau}=2,3,4$ can be seen in Table \ref{ButcherTableau}.

\begin{table}[htb]
\renewcommand\arraystretch{1.2}
\centering
\begin{subtable}[t]{0.3\textwidth}
\[
\begin{array}
{c|cccc}
0 & \frac{1}{2} & -\frac{1}{2}\\
1 & \frac{1}{2} & \frac{1}{2}\\
\hline
& \frac{1}{2} &\frac{1}{2} 
\end{array}
\]
\caption{$N_{\tau}=2$}
\end{subtable}
    \hfil
\begin{subtable}[t]{0.3\textwidth}
\[
\begin{array}
{c|cccc}
0&\frac{1}{6} &-\frac{1}{3} &\frac{1}{6} \\
\frac{1}{2} & \frac{1}{6} &\frac{5}{12} &-\frac{1}{12} \\
1& \frac{1}{6} &\frac{2}{3} &\frac{1}{6} \\
\hline
& \frac{1}{6} &\frac{2}{3} &\frac{1}{6} 
\end{array}
\]
\caption{$N_{\tau}=3$}
\end{subtable}
\medskip
\begin{subtable}[t]{0.4\textwidth}
\[
\begin{array}
{c|cccc}
0 & \frac{1}{12} &-\frac{\sqrt{5}}{12} &\frac{\sqrt{5}}{12} &-\frac{1}{12} \\
\frac{1}{2}-\frac{\sqrt{5}}{10} & \frac{1}{12} &\frac{1}{4} &\frac{10-7\sqrt{5}}{60} &\frac{\sqrt{5}}{60} \\
\frac{1}{2}+\frac{\sqrt{5}}{10} &\frac{1}{12} &\frac{10+7\sqrt{5}}{60}&\frac{1}{4} &-\frac{\sqrt{5}}{60} \\
1& \frac{1}{12} &\frac{5}{12} &\frac{5}{12} &\frac{1}{12} \\
\hline
& \frac{1}{12} &\frac{5}{12} &\frac{5}{12} &\frac{1}{12} 
\end{array}
\]
\caption{$N_{\tau}=4$}
\end{subtable}
\caption{Butcher Tableaus for Lobatto IIIC methods}\label{ButcherTableau}
\end{table}


\section*{Appendix C}\label{Appendix_UFL}
This patch adds 4D simplex and cuboid reference elements to UFL needed for the
3D$+ t$ simulations. This patch is currently implemented in \dunefem and will
be discussed with the UFL community. \\

\begin{python}
# 4d patching of reference elements 
def _patchufl4d():
    from ufl.sobolevspace import H1
    from ufl.finiteelement.elementlist import ufl_elements, any_cell, register_element
    from ufl.cell import num_cell_entities, cellname2facetname, 
    from ufl.cell import _simplex_dim2cellname, _hypercube_dim2cellname

    # check if this has been added before
    if not 'pentatope' in ufl.cell.num_cell_entities:
        # 4d-simplex
        ufl.cell.num_cell_entities["pentatope"] = (5, 10, 10, 5, 1) 
        # 4d-cube
        ufl.cell.num_cell_entities["tesseract"]   = (16, 32, 24, 8, 1)

        # recompute cell name to dimension mapping
        ufl.cell.cellname2dim = dict((k, len(v) - 1) for k, v \ 
          in ufl.cell.num_cell_entities.items())

        ufl.cell.cellname2facetname["pentatope"] = "tetrahedron"
        ufl.cell.cellname2facetname["tesseract"]   = "hexahedron"

        ufl.cell._simplex_dim2cellname[4]   = "pentatope"
        ufl.cell._hypercube_dim2cellname[4] = "tesseract"

        # add types to element lists
        ufl.finiteelement.elementlist.simplices =\
            ufl.finiteelement.elementlist.simplices + ("pentatope",)
        ufl.finiteelement.elementlist.cubes = \
            ufl.finiteelement.elementlist.cubes + ("tesseract",)
        ufl.finiteelement.elementlist.any_cell =\
                ufl.finiteelement.elementlist.any_cell + ("pentatope", "tesseract", )

        # register Lagrange again with new element type list
        ufl_elements.pop("Lagrange")
        ufl_elements.pop("CG")
        register_element("Lagrange", "CG", 0, H1, "identity", (1, None), \
          ufl.finiteelement.elementlist.any_cell)  
\end{python}

\begin{python}
# selecting a cell based on the dimension of the domain and or grid
def cell(dimDomainOrGrid):
    if isinstance(dimDomainOrGrid,ufl.Cell):
        return dimDomainOrGrid
    try:
        dimWorld = int(dimDomainOrGrid.dimWorld)
        dimDomain = int(dimDomainOrGrid.dimGrid)
    except:
        dimDomain = dimDomainOrGrid
        if isinstance(dimDomain, tuple):
            if len(dimDomain) != 2:
                raise Exception('dimDomain tuple must contain exactly two elements.')
            dimWorld = int(dimDomain[1])
            dimDomain = dimDomain[0]
        else:
            dimWorld = int(dimDomain)
    if dimDomain == 1:
        return ufl.Cell("interval", dimWorld)
    elif dimDomain == 2:
        return ufl.Cell("triangle", dimWorld)
    elif dimDomain == 3:
        return ufl.Cell("tetrahedron", dimWorld)
    elif dimDomain == 4:
        # add 4d cell types to ufl data structures
        _patchufl4d()
        return ufl.Cell("pentatope", dimWorld)
    else:
        raise NotImplementedError('UFL cell not implemented for dimension '\
           + str(dimDomain) + '.')
\end{python}


\section*{Appendix D}\label{Appendix_DGSEM}
For $N_{\tau}=2$ we get
\begin{align}
\mat{B}_{\tau} = 
\begin{pmatrix}
-1 & 0 \\
0 & 1
\end{pmatrix},~
\mat{M}_{\tau} = 
\begin{pmatrix}
1 & 0 \\
0 & 1
\end{pmatrix},~
\mat{D}_{\tau} = 
\begin{pmatrix}
-\frac{1}{2} & \frac{1}{2} \\
-\frac{1}{2} & \frac{1}{2}
\end{pmatrix}.
\end{align}
For $N_{\tau}=3$ we get
\begin{align}
\mat{B}_{\tau} = 
\begin{pmatrix}
-1 & 0 & 0\\
0 & 0 & 0 \\
0 & 0 & 1
\end{pmatrix},~
\mathbf{M}_{\tau} = 
\begin{pmatrix}
\frac{1}{3} & 0 & 0 \\
0 & \frac{4}{3} & 0 \\
0 & 0 & \frac{1}{3}
\end{pmatrix},~
\mat{D}_{\tau} = 
\begin{pmatrix}
-\frac{3}{2} & 2 & -\frac{1}{2} \\
-\frac{1}{2} & 0 & \frac{1}{2} \\
\frac{1}{2} & -2 & \frac{3}{2}
\end{pmatrix}.
\end{align}
For $N_{\tau}=4$ we get
\begin{align}
\mat{B}_{\tau} = 
\begin{pmatrix}
-1 & 0 & 0 & 0\\
0 & 0 & 0 & 0 \\
0 & 0 & 0 & 0 \\
0 & 0 & 0 & 1
\end{pmatrix},~
\mat{M}_{\tau} = 
\begin{pmatrix}
\frac{1}{6} & 0 & 0 & 0 \\
0 & \frac{5}{6} & 0 & 0 \\
0 & 0 & \frac{5}{6} & 0 \\
0 & 0 & 0 & \frac{1}{6}
\end{pmatrix},~
\mat{D}_{\tau} = 
\begin{pmatrix}
-3 & \frac{5+5\sqrt{5}}{4} &  \frac{5-5\sqrt{5}}{4} & \frac{1}{2}\\
\frac{-1-\sqrt{5}}{4} & 0 & \frac{-\sqrt{5}}{2} & \frac{1-\sqrt{5}}{4} \\
\frac{-1+\sqrt{5}}{4} & \frac{-\sqrt{5}}{2} & 0 & \frac{1+\sqrt{5}}{4} \\
-\frac{1}{2} & \frac{5\sqrt{5}-5}{4} & \frac{-5-5\sqrt{5}}{4} & 3
\end{pmatrix}.
\end{align}


\bibliographystyle{abbrv}
\bibliography{Sections/References.bib}


\end{document}